\documentclass[10pt]{article}
\usepackage{amsfonts}
\usepackage[leqno]{amsmath}
\usepackage{graphicx}
\usepackage{latexsym}
\usepackage{amsmath,amsfonts,amssymb,amsthm,mathrsfs,euscript,makeidx,color}
\usepackage{enumerate}
\allowdisplaybreaks[1]
\usepackage{multirow}

\def\sqr#1#2{{\vcenter{\vbox{\hrule height.#2pt
				\hbox{\vrule width.#2pt height#1pt \kern#1pt \vrule width.#2pt}
				\hrule height.#2pt}}}}
\def\signed #1{{\unskip\nobreak\hfil\penalty50
		\hskip2em\hbox{}\nobreak\hfil#1
		\parfillskip=0pt \finalhyphendemerits=0 \par}}
\def\endpf{\signed {$\sqr69$}}

\def\3n{\negthinspace \negthinspace \negthinspace }
\def\2n{\negthinspace \negthinspace }
\def\1n{\negthinspace }
\def\bel{\begin{equation}\label}

\def\dbE{\mathbb{E}}
\def\dbF{\mathbb{F}}

\def\dbP{\mathbb{P}}

\def\dbR{\mathbb{R}}
\def\dbS{\mathbb{S}}

\def\sK{\mathscr{K}}

\def\sV{\mathscr{V}}

\def\ds{\displaystyle}

\def\ns{\noalign{\ss}}
\def\a{\alpha}
\def\b{\beta }

\def\d{\delta}
\def\e{\varepsilon}
\def\z{\zeta}
\def\k{\kappa}
\def\l{\lambda}
\def\m{\mu}

\def\si{\sigma}
\def\t{\tau}
\def\f{\varphi}
\def\th{\theta}
\def\o{\omega}

\def\i{\infty}
\def\D{\Delta}

\def\Si{\Sigma}
\def\F{\Phi}
\def\O{\Omega}
\def\cF{{\cal F}}

\def\BH{{\bf H}}

\def\D{\Delta}

\def\Si{\Sigma}
\def\F{\varPhi}
\def\O{\Omega}

\def\BBXi{\boldsymbol\Xi}

\def\no{\noindent}

\def\ss{\smallskip}
\def\ms{\medskip}

\def\q{\quad}
\def\qq{\qquad}
\def\hb{\hbox}

\def\h1{\outline{$1$}}

\def\hh1{\outline{$1$}}
\def\hh2{\outline{$2$}}
\def\hh3{\outline{$3$}}
\def\hh4{\outline{$4$}}
\def\hh5{\outline{$5$}}
\def\hh6{\outline{$6$}}
\def\hh7{\outline{$7$}}
\def\hh8{\outline{$8$}}
\def\hh9{\outline{$9$}}
\def\hh0{\outline{$0$}}

\def\lan{{\langle}}
\def\ran{{\rangle}}

\def\h{\widehat}
\def\wt{\widetilde}

\def\cd{\cdot}
\def\cds{\cdots}

\def\as{\hbox{\rm a.s.}}

\def\tr{\hbox{\rm tr$\,$}}

\def\les{\leqslant}
\def\ges{\geqslant}

\def\({\Big (}
\def\){\Big )}
\def\[{\Big[}
\def\]{\Big]}
\def\lan{\langle}
\def\ran{\rangle}

\def\bde{\begin{definition}\label}
	\def\ede{\end{definition}}
\def\bel{\begin{equation}\label}
		\def\ee{\end{equation}}
	\def\bt{\begin{theorem}\label}
		\def\et{\end{theorem}}
	\def\bc{\begin{corollary}\label}
		\def\ec{\end{corollary}}
	\def\bl{\begin{lemma}\label}
		\def\el{\end{lemma}}
	\def\bp{\begin{proposition}\label}
		\def\ep{\end{proposition}}
	\def\bex{\begin{example}\label}
		\def\ex{\end{example}}
	\def\bas{\begin{assumption}}
		\def\eas{\end{assumption}}
	\def\br{\begin{remark}\label}
		\def\er{\end{remark}}
	\def\ba{\begin{array}}
		\def\ea{\end{array}}
	\def\ed{\end{document}}

\def\rf{\eqref}

\def\square#1{\vbox{\hrule\hbox{\vrule height#1%
			\kern#1\vrule}\hrule}}
\def\rectangle#1#2{\vbox{\hrule\hbox{\vrule height#1%
			\kern#2\vrule}\hrule}}

\font\tenbb=msbm10 \font\sevenbb=msbm7 \font\fivebb=msbm5

\newfam\bbfam
\scriptscriptfont\bbfam=\fivebb \textfont\bbfam=\tenbb
\scriptfont\bbfam=\sevenbb

\newtheorem{theorem}{Theorem}[section]
\newtheorem{corollary}[theorem]{Corollary}

\newtheorem{lemma}[theorem]{Lemma}
\newtheorem{proposition}[theorem]{Proposition}

\theoremstyle{definition}
\newtheorem{definition}[theorem]{Definition}

\newtheorem{remark}[theorem]{Remark}

\newtheorem{example}[theorem]{Example}

\makeatletter

\@addtoreset{equation}{section}
\makeatother

\usepackage{xcolor}
\makeatletter

\input pdf-trans
\newbox\qbox
\def\usecolor#1{\csname\string\color@#1\endcsname\space}
\newcommand\bordercolor[1]{\colsplit{1}{#1}}
\newcommand\fillcolor[1]{\colsplit{0}{#1}}
\newcommand\outline[1]{\leavevmode%
	\def\maltext{#1}%
	\setbox\qbox=\hbox{\maltext}%
	\boxgs{Q q 2 Tr \thickness\space w \fillcol\space \bordercol\space}{}%
	\copy\qbox%
}
\makeatother
\newcommand\colsplit[2]{\colorlet{tmpcolor}{#2}\edef\tmp{\usecolor{tmpcolor}}%
	\def\tmpB{}\expandafter\colsplithelp\tmp\relax%
	\ifnum0=#1\relax\edef\fillcol{\tmpB}\else\edef\bordercol{\tmpC}\fi}
\def\colsplithelp#1#2 #3\relax{%
	\edef\tmpB{\tmpB#1#2 }%
	\ifnum `#1>`9\relax\def\tmpC{#3}\else\colsplithelp#3\relax\fi
}
\bordercolor{black}
\fillcolor{white}
\def\thickness{.3}

\begin{document}

\title{\bf Stochastic Optimal Impulse Controls with\\ Changing Running Costs}

\author{Yuchen Cao\thanks{Citibank, N.A., Irving, TX, 75063, USA
                            (Email: {\tt yuchen1.cao@citi.com}).
                           }
~~~and~~~
Jiongmin Yong\thanks{Department of Mathematics, University of Central Florida, Orlando, FL 32816, USA
                    (Email: {\tt jiongmin.yong@ucf.edu}).
                    This author is supported by NSF grant DMS-2305475.}
}

\date{}
\maketitle

\ms

\no{\bf Abstract.} This paper is concerned with stochastic impulse control problems in which the running cost changes depending on the impulse control. Because of such a dependence, it brings several difficulties when the usual dynamic programming principle is to be used. The corresponding  Hamilton-Jacobi-Bellman (HJB) equation (a quasi-variational inequality) is derived, which contains a parameter. The value function is a unique viscosity solution to this HJB equation by a classical argument. Further, inspired by the derivation of the Pontryagin type maximum principle for stochastic optimal controls with a non-convex control domain, we have established the maximum principle for our stochastic optimal impulse controls, allowing perturbations in optimal impulse moments.

\ms

\no{\bf Keywords.} Optimal impulse control, dynamic programming principle, Hamilton-Jacobi-Bellman equations, quasi-variational inequality, viscosity solutions, maximum principle.

\ms

\no{\bf AMS 2020 Mathematics Subject Classification.}  93E20, 93C27, 49L25, 49N25.

\section{Introduction.}\label{Sec:Intro}

Let $(\O,\cF,\dbF,\dbP)$ be a filtered complete probability space satisfying the usual conditions. In this probability space, we define a $d$-dimensional standard Brownian motion $W(\cd)$ with $\dbF$ being its natural filtration, adding all the $\dbP$-null sets.

\ms

We begin with a motivation. Let $r(\cd)$ be the interest rate process modeled by a certain SDE, say, Cox-Ingersoll-Ross (CIR, for short) see \cite{Altay-Schmock 2013}. Let $t\in[0,T)$ be the current time. Suppose a unit (a company, for example) is to borrow certain amount of money $\xi_j$ at different moments $\t_j\ges t$ at (possibly) different rate $r(\t_j)>0$ from a bank. We assume the unit is getting into the system at $t$. The hard part is that $t$ might not be a moment at which a loan amount is initiated. But at that moment, a loan amount might be borrowed (at a rate settled at some moment $\t_0<t$). Now, since at $\t_j\ges t$, a moment at which a loan amount $\xi_j$ is initiated, the unit pays the cost $\ell(\t_j,\xi_j)>0$. Then the total amount owed by the unit is
$$\ba{ll}
\ns\ds C(s)=e^{r(\t_0)(s-t)}\xi_0+\ell(\t_1,\xi_1)+e^{r(\t_1)(s-\t_1)}\xi_1
+\cds\\
\ns\ds\qq=e^{r(\t_0)(s-t)}\xi_0+\sum_{j\ges1}\(\ell(\t_j,\xi_j)+e^{r(\t_j)
(s-\t_j)}\xi_j\){\bf1}_{[\t_j,T]}(s),\qq s\in[t,T].\ea$$
Or equivalently,
$$C(s)=\xi_0+\int_t^sr(\t_0)\xi_0d\t+\sum_{j\ges1}\(\ell(\t_j,\xi_j)+\xi_j
+\int_{\t_j}^sr(\t_j)\xi_jd\t\){\bf1}_{[\t_j,T]}(s).$$
Now, with these loans, the unit has some investments to get better returns. So it is reasonable to assume the change in the wealth of the unit satisfies the SDE:
\bel{state}\left\{\2n\ba{ll}
\ds dX(s)=\big[b(s,X(s))+\xi(s)\big]ds+\si(s,X(s))dW(s),\qq s\in[t,T],\\
\ns\ds X(t)=x,\ea\right.\ee
with
\bel{xi}\xi(s)=\sum_{j\ges1}\xi_j{\bf1}_{[\t_j,T]}(s),\qq s\in[t,T].\ee
Then the expected utility (profit) on $[t,T]$ (of the unit) is
\bel{cost}J(t,x;\xi(\cd))=\dbE\[\int_t^Tg(\t_0,\th,X(\th))d\th
+\sum_{j\ges1}\(\int_{\t_j}^Tg(\t_j,\th,X(\th))d\th+\ell(\t_j,
\xi_j)\)+h(X(T))\].\ee
The dependence on all the $\t_j$ is to be emphasized because on each interval $[\t_j,T]$, the payment of the unit to the bank should be different, and the better return might be different as well. We may formulate an optimal control problem with the state equation \rf{state} and payoff functional \rf{cost}, with respect to the control $\xi(\cd)$.

\ms

Let us present another motivation. Still, let $t\in[0,T)$ be the current time. Suppose a unit (a company, for example, again) needs to buy several raw materials (which could be contained in some other material, for example, some iron-containing minerals such as hematite, magnetite, goethite. limonite, and siderite, etc.). Let $\xi_j$ be the raw material (which could be some organic or inorganic ones) bought at moment $\t_j\ges t$, with cost $\ell(\t_j,\xi_j)$. Then the unit owns the raw material of form \rf{xi} with each $\xi_j$ possibly being a vector (of the same dimensions as the state). The unit produces certain products using the raw material $\xi(s)$. Thus, the state process $X(\cd)$, satisfying the state equation \rf{state}, represents the quantity of the products for the unit. The expected payoff functional is of the form \rf{cost}. Here, the reason for $g$ depending on $\t_j$ is that by having $\xi_j$ at $\t_j$ (from other raw material which could get $\xi_j$ differently), the running payoff will be different.

\ms

The change in running cost also happens when people consider mortgage refinance related problems. Note that main motivations for refinancing are the following: the interest rate decreases significantly, the monthly payment is largely reduced, some cash can be further borrowed (according to the credit line), the pre-payment of the mortgage can be accelerated, etc. In any of these cases, the new mortgage will follow the new mortgage rate. Thus, it is natural that the running cost rate will be changed after the new mortgage starts.

\ms

Based on the above motivations, we now formulate our optimal impulse control problem in general.

\ms

Let $\sK\subseteq\dbR^n$ be a closed and convex cone, and let
$$\D_*[0,T]=\Big\{(\t_0,t)\in[0,T]\times[0,T]\bigm|\t_0\les t\Big\}$$
be the lower triangle domain. Let $t\in[0,T)$ be given (thus, $t$ is the moment that we get into the system) and $\BBXi[t,T]$ be the set of all impulse controls of the form \rf{xi} with $\xi_j\in\sK$ (which might be a multi-dimensional vector, since $n$ might be larger than 1), $\cF_{\t_j}$-measurable and integrable. For any $\xi(\cd)\in\BBXi[t,T]$, the state equation looks like the following:
\bel{state}\ba{ll}
\ns\ds X(s)=x+\sum_{j\ges0}\(\int_t^s\3n b(\t_j,\th,X(\th)){\bf1}_{[\t_j,T]}(\th)d\th
+\2n\int_t^s\3n\si(\t_j,\th,X(\th)){\bf1}_{[\t_j,T]}(\th)dW(\th)\)\1n+\1n\xi(s),\q s\in[t,T],\ea\ee
where $b:\D_*[0,T]\times\dbR^n\to\dbR^n$ and $\si:\D_*[0,T]\times\dbR^n\to
\dbR^{n\times d}$ are given maps. The cost functional is
\bel{J}J(\t_0,t,x;\xi(\cd))=\dbE\[\bar Y(t,\t_0,t,x;\xi(\cd))\],\ee
where $(\bar Y(\cd),\bar Z(\cd\,,\cd))\equiv(\bar Y(\cd\,;\t_0,t,x,\xi(\cd)),\bar Z(\cd\,,\cd\,;
\t_0,t,x,\xi(\cd))$ is the adapted solution of the following backward stochastic
Volterra integral equation (BSVIE, for short):
\bel{BSVIE}\ba{ll}
\ns\ds\bar Y(s)=h(X(T))+\sum_{j\ges1}\ell(\t_j,\xi_j)+\int_s^T\sum_{j\ges0}
g(\t_j,\th,X(\th)){\bf1}_{[\t_j,T]}(\th)d\th\\
\ns\ds\qq\qq\qq\qq-\int_s^T\bar Z(s,\th)dW(\th),\qq s\in[t,T].\ea\ee
Note that in the above, the free term
$$\psi(s)\equiv h(X(T))+\sum_{\t_j\ges s}\ell(\t_j,\xi_j),\qq s\in[t,T]$$
is an $\cF_T$-measurable process (not just an $\cF_T$-measurable random variable) because of the second summation term. Clearly, by taking $s=t$ and applying the operator $\dbE$ on both sides of \rf{BSVIE}, we get
\bel{J*}J(\t_0,t,x;\xi(\cd))=\dbE\[h(X(T))+\sum_{\t_j\ges s}\ell(\t_j,\xi_j)
+\int_t^T\sum_{j\ges1}g(\t_j,\th,X(\th)){\bf1}_{[\t_j,T]}(\th)d\th\].\ee
This shows that our cost functional is an extension of classical one. Our optimal impulse control problem can be formulated as follows.

\ms

{\bf Problem (IC).} For given $(\t_0,t,x)\in\D_*[0,T]\times\sK$, find $\bar\xi(\cd)\in\BBXi[t,T]$ such that
\bel{J=V}J(\t_0,t,x;\bar\xi(\cd))=\inf_{\xi(\cd)\in\BBXi[t,T]}J(\t_0,t,x;\xi(\cd))
\equiv V(\t_0,t,x).\ee

\ms

We see that our stochastic optimal impulse control problem is different from the classical one in the following aspects:

\ms

$\bullet$ The drift and the diffusion terms in the state equation \rf{state} depend on $\t_0\les t$, which is no later than the current moment $t$ (and might be before the current moment). This brings an obvious difficulty for people to use dynamic programming principle (DPP, for short) because people might get into the system at time which could not be an impulse moment.

\ms

$\bullet$ The cost functional \rf{J} is recursive involving a BSVIE, which contains the impulse costs in the free term which might be an $\cF_T$-measurable process (and not necessarily a $\cF_T$-measurable random variable or a $\dbF$-adapted process).

\ms

It is well known that impulse control problems were initiated in 1973 by Bensoussan--Lions (\cite{Bensoussan-Lions 1973}). According to Remark 3 of \cite{Bensoussan-Lions 1973}, they reformulated the so-called inventory problem in an earlier work of Bather \cite{Bather 1966} and formally introduced quasi-variational inequalities so that a new research area was opened.\footnote{The authors would like thank Professor Hanxiao Wang of Shenzhen University (China) for his help to clarify this fact.} There are many other follow up works. We can only list a small portion here. For determined version, see \cite{Barles 1985, Rempala-Zabczyk 1988, Yong-Zhang 1992, Silva-Vinter 1997, El Farouq-Barles-Berhnard 2010}, For SDE systems, relevant problems, and applications, see \cite{Bensoussan-Frehse-Mosco 1982, Perthame 1985, Frehse-Mosco 1982, Matzeu-Mosco-Vivaldi 1983, Bensoussan-Lions 1984, Eastham-Hasting 1988, Tang-Yong 1993, Korn 1997, Korn 1999, Yang 2001, Ohnishi-Tsujimura 2006, Egami 2008, Oksendal-Sulem 2008, Alarez-Lempa 2008, Bruder-Pham 2009, Djehuche-Hamadene-Hdhiri 2010, Leander-Lenhart-Protopopescu 2015, Jain 2016, Menaldi-Robin 2017, Hdhiri-Karouf 2018, Piunovskiy-et-al 2019, Belak-Christensen-Seifried 2017, Feng-Muthuraman 2010, El Asri-Mazid 2020, Kurzhanski-Daryin 2020, Li-Yong 2021, Djehuche-Hamadene-Hdhiri-Zaatra 2022, Lv-Xiong 2022, Josson-Perninge 2023, Fusco-Motta-Vinter 2024, Perninge 2024}. For relevant differential games, see \cite{Yong 1994, Zhang 2011, Perninge 2023}. We note that the above stochastic impulse control problems were treated by dynamic programming principle and PDE methods. As we know that one can approach the problems by variational method as well. Along this line, we mention the literature \cite{Rempala-Zabczyk 1988, Yong-Zhang 1992, Silva-Vinter 1997} for deterministic versions of the maximum principle, and \cite{Hu-Yong 1991, Yu 2012, Wu-Zhang 2012, Wu-Zhang 2022, Wu-Zhang 2024, Yang-Yu 2025} for the stochastic counterpart. See below for more comments.

\ms

Due to the existence of $\t_0$ (see \rf{J=V}), we still have the continuity of the value function $V(\t_0,t,x)$ by making use of the monotonicity of the impulse cost. The HJB equation for the value function contains $\t_0$ as a parameter. Then, the classical definition of viscosity solution for such a HJB equation has to be modified. Under such a definition, the classical argument is still applicable and the value function is the unique viscosity solution to this HJB equation. This gives the first contribution of the current paper.

\ms

For optimal (impulse) controls, as we indicated earlier that one may use variational technique to derive Pontryagin type maximum principle. In \cite{Hu-Yong 1991}, a classical stochastic optimal impulse control problem (with fixed running cost) was considered. To derive the maximum principle, in order to get the right order of the variational system, only the sizes of the optimal impulse were perturbed (See \cite{Yu 2012, Wu-Zhang 2012, Wu-Zhang 2022, Wu-Zhang 2024, Yang-Yu 2025} also). It turns out that to perturb the optimal impulse moments, the spike variational technique introduced by Peng \cite{Peng 1990} (See also \cite{Yong-Zhou 1999, Yong 2022}) is needed, and the second order adjoint equation has to be used. This is another main contribution of the current paper.

\ms

The rest of the paper is organized as follows. In Section 2, some properties of the state process and the cost functional are presented. Dynamic programming principle for the value function and related results are collected in Section 3. The maximum principle of Pontryagin type for our optimal impulse control has been established in Section 4. Some possible extensions are indicated in Section 5. Several concluding remarks are put in that section as well.

\section{Properties of State Processes and Cost Functionals.}

We first make some observations. For any $(t,x)\in[0,T]\times\dbR^n$, let
$\xi(\cd)\in\BBXi[t,T]$ of form \rf{xi}. Denote the number of impulses in $\xi(\cd)$ by $\k(\xi(\cd))$. Due to the positivity of the impulse cost (we are going to assume), we must have that $\k(\xi(\cd))<\i$. Let $\t_0\in[0,t]$. Then from \rf{state}, with $j_0\ges1$, we have
\bel{X(t)}\ba{ll}
\ns\ds X(\t_{j_0})=X(\t_{j_0}+0)=x+\sum_{j\ges0}\(\int_t^{\t_{j_0}}b(\t_j,\th,X(\th))
{\bf1}_{[\t_j,T]}(\th)d\th\\
\ns\ds\qq\qq\qq\qq\qq\qq\qq+\int_t^{\t_{j_0}}\si(\t_j,\th,X(\th))
{\bf1}_{[\t_j,T]}(\th)dW(\th)\)+\sum_{j=1}^{j_0}\xi_j{\bf1}_{[\t_j,T]}(s),\ea\ee
Hence, equivalent to \rf{state}, for $\t_{j_0}\les s<\t_{j_0+1}$, we have
\bel{state*} X(s)=X(\t_{j_0})+\sum_{j=0}^{j_0}\(\int_{\t_{j_0}}^sb(\t_j,\th,X(\th))d\th
+\int_{\t_{j_0}}^s\si(\t_j,\th,X(\th))dW(\th)\).\ee
Note that we are going to prove the fact that in the interval $(\t_{j_0},\t_{j_0+1})$, there is no impulse, and the impulse is contained in the term $X(\t_{j_0})$. For $s\in[t,\t_1)$,
\bel{state**}X(s)=x+\int_t^sb(\t_0,\th,X(\th))d\th
+\int_t^s\si(\t_0,\th,X(\th))dW(\th).\ee

Now, we introduce the following assumption.

\ms

{\bf(H1)} Let $b,\si:\D_*[0,T]\times\dbR^n\to\dbR^n$
be continuous, for some generic constant $K>0$, and some generic modulus of continuity $\o:\dbR_+\to\dbR_+$,
\bel{b}|b(\t,t,0)|+|\si(\t,t,0)|\les K,\qq(\t,t)\in\D_*[0,T].\ee
\bel{b-b}\ba{ll}
\ns\ds|b(\t,t,x)-b(\t',t,x')|+|\si(\t,t,x)-\si(\t',t,x')|\les\o(|\t-\t'|)+ K|x-x'|,\\
\ns\ds\qq\qq\qq\qq\qq\qq\qq\qq(\t,t,x),(\t',t,x')\in\D_*[0,T]\times\dbR^n.\ea\ee

\ms

In what follows, $K>0$ will represent a generic constant and $\o(\cd)$ will represent a generic modulus of continuity, which could be different from line to line, unless otherwise stated. For convenience, we let $d=1$ (general case can be treated in a similar way) and some conditions in (H1) can be relaxed (we assume  those just for convenience). Also, for any $\xi(\cd)\in\BBXi[t,T]$ of form \rf{xi}, we denote the number of impulses by $\k(\xi(\cd))$. The following result is about the controlled process.

\bp{prop2.1} \sl Let {\rm(H1)} hold. Then for any $(t,x)\in[0,T]\times\dbR^n$, $\xi(\cd)\in\BBXi[t,T]$, and $\t_0\in[0,t]$, state equation \rf{state}
admits a unique solution $X(\cd)=X(\cd\,;\t_0,t,x,\xi(\cd))$.
Moreover, for $p\ges1$, there exists a constant $K>0$, depending on the number $\k(\xi(\cd))$ of the impulses in $\xi(\cd)$ such that
\bel{X}\ba{ll}
\ns\ds\dbE\[\sup_{s\in[t,T]}|X(s;\t_0,t,x,\xi(\cd))|^p\]\\
\ns\ds\les K\dbE\Big\{\sum_{j=1}^{\k(\xi(\cd))}|\xi_j|^p+\sum_{j=0}^{\k(\xi(\cd))}\[\(\int_{\t_j}^{\t_{j+1}}|b(\t_j,\th,0)|d\th\)^p
\2n+\(\int_{\t_j}^{\t_{j+1}}\2n|\si(\t_j,\th,0)|^2d\th\)^{p\over2}\]\\
\ns\ds\qq\qq+|x|^p+\(\int_t^{\t_1}|b(\bar\t_0,\th,0)|d\th\)^p+\(\int_t^{\t_1}
|\si(\bar\t_0,\th,0)|^2d\th\)^{p\over2}\\
\ns\ds\qq\qq+\(\int_{t_{\k(\xi(\cd))}}^T|b(\t_{\k(\xi(\cd))},\th,0)|d\th\)^p
+\(\int_{\t_{\k(\xi(\cd))}}^T|\si(\t_{\k(\xi(\cd))},\th,0)|^2
d\th\)^{p\over2}\Big\}.\ea\ee
where $\bar\t_0=\t_0$. If $(\t_0',x')\in[0,t]\times\dbR^n$, then for some $K>0$, depending on $\k(\xi(\cd))$,
\bel{X-X}\dbE\[\sup_{s\in[t,T]}|X(s;\t_0,t,x,\xi(\cd))-X(s;\t_0',t,x',\xi(\cd))|^p\]
\les\o(|\t_0-\t_0'|)+K|x-x'|^p.\ee

\ep

\it Proof. \rm First of all, for any $(t,x)\in[0,T]\times\dbR^n$ and $\xi(\cd)\in\BBXi[t,T]$ (of
form \rf{xi}), by recursive procedure, we have the unique strong solution to the state equation \rf{state}. Further, for $\t_{j_0}\les s<\t_{j_0+1}$, with $j_0\ges1$, we have (from \rf{state*})
$$\ba{ll}
\ns\ds\dbE\[\sup_{r\in[\t_{j_0},s]}|X(r)|^p\]\les\dbE
\[\sup_{r\in[\t_{j_0},s]}\Big|X(\t_{j_0})+\sum_{j=0}^{j_0}
\(\int_{\t_{j_0}}^r\3n b(\t_j,\th,X(\th))d\th+\int_{\t_{j_0}}^r\3n\si(\t_j,\th,X(\th))
dW(\th)\)\Big|^p\]\\
\ns\ds\les(2j_0+1)^{p-1}\dbE\Big\{|X(\t_{j_0})|^p+\sum_{j=0}^{j_0}
\[\(\int_{\t_{j_0}}^s\3n|b(\t_j,\th,X(\th))|d\th\)^p+\sup_{r\in[\t_{j_0},s]}
\Big|\int_{\t_{j_0}}^r\3n\si(\t_j,\th,X(\th))dW(\th)\Big|^p\]\Big\}\\
\ns\ds\les(2j_0+1)^{p-1}\dbE\Big\{|X(\t_{j_0})|^p+\sum_{j=0}^{j_0}
\[\(\int_{\t_{j_0}}^s\3n|b(\t_j,\th,X(\th))|d\th\)^p+K_p\(\int_{\t_{j_0}}^s
|\si(\t_j,\th,X(\th))|^2d\th\)^{p\over2}\Big|^p\]\Big\}\\
\ns\ds\les K\dbE\Big\{|X(\t_{j_0})|^p+\sum_{j=0}^{j_0}\[\(\int_{\t_{j_0}}^s
\big[|b(\t_j,\th,0)|+K|X(\th)|\big]d\th\)^p+\(\int_{\t_{j_0}}^s
\big[|\si(\t_j,\th,0))|+K|X(\th)|\big]^2d\th\)^{p\over2}\]\Big\}\\
\ns\ds\les K\dbE\Big\{|X(\t_{j_0})|^p+\sum_{j=0}^{j_0}\[\(\int_{\t_{j_0}}^s|b(\t_j,\th,0)|d\th\)^p
\2n+\(\int_{\t_{j_0}}^s\2n|\si(\t_j,\th,0)|^2d\th\)^{p\over2}\]\Big\}+K\2n
\int_{\t_{j_0}}^s\2n\dbE\[\sup_{r\in[\t_{j_0},\th]}|X(r)|^p\]d\th.\ea$$
In the above, $K_p>0$ is the absolute constant appearing in the Burkholder-Davis-Gundy inequalities (see \cite{Karatzas-Shreve 1998}). Consequently, by Gronwall's inequality, we have
$$\dbE\[\sup_{r\in[\t_{j_0},\t_{j_0+1})}|X(r)|^p\]\les K\dbE\Big\{|X(\t_{j_0})|^p+\sum_{j=0}^{j_0}\[\(\int_{\t_{j_0}}^s|b(\t_j,\th,0)|d\th\)^p
\2n+\(\int_{\t_{j_0}}^s\2n|\si(\t_j,\th,0)|^2d\th\)^{p\over2}\]\Big\},\q j_0\ges1,$$
with $K$ depending on $\k(\xi(\cd))$. Note
$$X(\t_{j_0})=X(\t_{j_0}-0)+\xi_{j_0},\q\hb{ or }\q\D X(\t_{j_0})=X(\t_{j_0})-X(\t_{j_0}-0)=\xi_{j_0}.$$
Also, for $r\in[t,\t_1)$,
$$\dbE\[\sup_{r\in[t,\t_1)}|X(r)|^p\]\les K\[|x|^p+\(\int_t^{\t_1}|b(\t_0,\th,0)|d\th\)^p+\(\int_t^{\t_1}|\si(\t_0,\th,0)|^2
d\th\)^{p\over2}\].$$
For $s\in[\t_{\k(\xi(\cd))},T]$,
$$\dbE\[\sup_{r\in[t_{\k(\xi(\cd))},T]}|X(r)|^p\]\les K\[|X(\k(\cd))|^p+\(\int_{t_{\k(\xi(\cd))}}^T|b(\t_0,\th,0)|d\th\)^p
+\(\int_{\t_{\k(\xi(\cd))}}^T|\si(\t_0,\th,0)|^2
d\th\)^{p\over2}\].$$
Thus,
$$\ba{ll}
\ns\ds\dbE\[\sup_{r\in[t,T]}|X(r)|^p\]\\
\ns\ds=\max\Big\{\dbE\[\sup_{r\in[t,\t_1)}|X(r)|^p\],
\dbE\[\sup_{r\in[\t_1,\t_2)}|X(r)|^p\],\cds,\dbE\[\sup_{r\in[\t_{\bar j_0-1},\t_{\bar j_0})}|X(r)|^p\],\dbE\[\sup_{r\in[\t_{\bar j_0},T]}|X(r)|^p\]\Big\}\\
\ns\ds\les K\dbE\Big\{\sum_{j=1}^{\k(\xi(\cd))}|\xi_j|^p+\sum_{j=0}^{\k(\xi(\cd))}\[\(\int_{\t_j}^{\t_{j+1}}|b(\t_j,\th,0)|d\th\)^p
\2n+\(\int_{\t_j}^{\t_{j+1}}\2n|\si(\t_j,\th,0)|^2d\th\)^{p\over2}\]\\
\ns\ds\qq\qq+|x|^p+\(\int_t^{\t_1}|b(\t_0,\th,0)|d\th\)^p+\(\int_t^{\t_1}|\si(\t_0,\th,0)|^2
d\th\)^{p\over2}\\
\ns\ds\qq\qq+\(\int_{t_{\k(\xi(\cd))}}^T|b(\t_{\k(\xi(\cd))},\th,0)|d\th\)^p
+\(\int_{\t_{\k(\xi(\cd))}}^T|\si(\t_{\k(\xi(\cd))},\th,0)|^2
d\th\)^{p\over2}\Big\}\ea$$
Then, \rf{X} easily follows.

\ms

Similarly, for $(\t_0,t,x),(\t_0',t,x')\in\D_*[0,T]\times\dbR^n$, if we denote $X'(\cd)=X(\cd\,;\t_0',t,x',\xi(\cd))$, then for $s\in[\t_{j_0},\t_{j_0+1})$, with $j_0\ges1$,
$$\ba{ll}
\ns\ds\dbE\[\sup_{r\in[\t_{j_0},s)}|X(r)-X'(r)|^p\]\\
\ns\ds\les (2j_0+1)^{p-1}\Big\{|X(\t_{j_0})-X'(\t_{j_0})|^p+(j_0+1)\(\int_{\t_{j_0}}^s
\big[\o(|\t_0-\t_0'|)+K|X(\th)-X'(\th)|\big]d\t\)^p\\
\ns\ds\qq\qq\qq\qq\qq\qq+(j_0+1)K_p\(\int_{\t_{j_0}}^s
K^2|X(\th)-X'(\th)|^2d\th\)^{p\over2}\Big\}.\ea$$
Hence,
$$\dbE\[\sup_{r\in[\t_{j_0},s)}|X(r)-X'(r)|^p\]\les\o(|\t_0-\t_0'|)+K |X(\t_{j_0})-X'(\t_{j_0})|^p,$$
with $\o(\cd)$ and $K>0$ depending on $\k(\xi(\cd))$.
Now, for $s\in[t,\t_1)$, we have (similar to the above)
$$\ba{ll}
\ns\ds\dbE\[\sup_{r\in[t,s)}|X(r)-X'(r)|^p\]\les K\dbE\[|x-x'|^p+\(\int_t^s\big[\o(|\t_0-\t_0'|)+K|X(\th)-X'(\th)|\big]d\th\)^p\\
\ns\ds\qq\qq\qq\qq\qq\qq\qq\q+\(\int_t^s\big[\o(|\t_0-\t_0'|)+K|X(\th)-X'(\th)|\big]^2
d\th\)^{p\over2}\].\ea$$
Then, we have \rf{X-X}. \endpf

\ms

We have seen that \rf{X-X} gives the continuous dependence of the state process on $(\t,x)$ with the impulse control to be fixed. Due to the jumps of the impulse control, we do not expect to have continuity of the state process with respect to the running time $s$, nor the variation of the impulse controls. This brings some essential difficulties in deriving the necessary conditions for optimal impulse controls. See Section 4 for details. On the other hand, the continuity of the value function $V(\t_0,t,x)$ in $t$ (The continuity in $(\t_0,x)$ is relatively simple) can be obtained by the monotonicity of the impulse cost $\ell(\t,\xi)$ in $\t$. This will be made precise below.

\ms
Now, for the cost functional, we introduce the following assumption.

\ms

{\bf(H2)} The maps $g:\D_*[0,T]\times\dbR^n\to\dbR_+$, $h:\dbR^n\to\dbR_+$ and
$\ell:[0,T]\times\sK\to\dbR_+$ are continuous. For some constant $K>0$ and uniform modulus of continuity $\o:\dbR_+\to\dbR_+$, it holds
\bel{g}|g(\t,t,x)|+|h(x)|\les K,\qq(\t,t,x)\in\D_*[0,T]\times\dbR^n,\ee
\bel{g-g}|g(\t,t,x)-g(\t',t,x')|\les\o(|\t-\t'|)+K|x-x'|,\q(\t,t,x),
(\t',t,x')\in\D_*[0,T]\times\dbR^n.\ee
The maps $\t\mapsto\ell(\t,\xi)$ is non-increasing, and for some $\ell_0>0$ and $\m\in(0,1]$,
\bel{ell}\ba{ll}
\ns\ds|\ell(\t,\xi)-\ell(\t',\xi|\les\o(|\t-\t'|),\qq \t,\t'\in[0,T],~\xi\in\sK,\\
\ns\ds\ell(\t,\xi)\ges\ell_0(1+|\xi|^\m),\qq\q(\t,\xi)\in[0,T]\times\sK,\\
\ns\ds\ell(\t,\xi+\xi')<\ell(\t,\xi)+\ell(\t,\xi'),\qq \t\in[0,T],~\xi,\xi'\in\sK.\ea\ee

\ms

Again, we may relax some conditions in (H2). We introduce (H2) only for convenience. The main result of this section is the following.

\bp{V} \sl Let {\rm(H1)--(H2)} hold. Then the value function $V(\t_0,t,x)$  is continuous and of linear growth. More precisely,
\bel{0<V}0\les V(\t_0,t,x)\les K,\ee
with $K>0$ being a generic constant, and
\bel{V-V}|V(\t_0,t,x)-V(\t_0',t,x')|\les \o(|\t_0-\t_0'|)+K(|t-t'|^{1\over2}+|x-x'|),\q(\t_0,t,x),(\t_0',t,x)\in\D_*[0,T]\times\dbR^n.\ee

\ep

\it Proof. \rm First of all, since under (H2), $J(t,x;\xi(\cd))$ is non-negative. So, we have the non-negativity of the value function $V(\cd\,,\cd\,,\cd)$. Next, under the trivial impulse control $\xi_0(\cd)$, the cost functional satisfies
$$J(\t_0,t,x;\xi_0(\cd))=\dbE\[\int_t^Tg(\t_0,s,X(s))ds
+h(X(T))\]\les K.$$
Hence, \rf{0<V} follows.

\ms

From \rf{0<V}, we see that in order $\xi(\cd)\in\BBXi[t,T]$ to be optimal, this $\xi(\cd)$ has to be in the following set:
\bel{BBXi}\BBXi_R[t,T]=\Big\{\xi(\cd)=\sum_{j\ges1}\xi_j{\bf1}_{[\t_j,T]}(\cd)\in
\BBXi[t,T]\bigm|\k(\xi(\cd))+\max_{j\ges1}|\xi_j|\les R\Big\},\ee
for $R>0$ large enough (might be depending on $|x|$). Thus, actually, without loss of generality, we may assume the above modulus of continuity $\o(\cd)$ and constant $K>0$ are independent of $\xi(\cd)\in\BBXi[t,T]$.
Now, from Proposition \ref{prop2.1}, we see that with $\t_0<\t_0'$ (Denote $X'(\cd)=X(\cd\,;\t_0',t,x',\xi(\cd))$)
$$\ba{ll}
\ns\ds|J(\t_0,t,x;\xi(\cd))-J(\t_0',t,x',\xi(\cd))|\les\dbE\[|h(X(T))-h(X'(T))|
+\int_t^T|g(\t_0,\th,X(\th))-g(\t_0',\th,X'(\th))|d\th\\
\ns\ds\qq\qq\qq\qq\qq\qq\qq\qq\qq+\int_t^T\sum_{j\ges1}|g(\t_j,\th,X(\th))-g(\t_j,\th,X'(\th))|d\th\]\\
\ns\ds\les\dbE\[K|X(T)-X'(T)|+\int_t^T\(\o(|\t_0-\t_0'|)+(\,\h{\1n j}+1)K|X(\th)-X'(\th)|\)d\th\]\les\o(|\t_0-\t_0'|)+K|x-x'|.\ea$$
Hence, \rf{V-V} follows with $t=t'$.

\ms

Finally, let $(\t_0,t,x),(\t_0,t',x)\in\D_*[0,T]\times\dbR^n$, with $t<t'$. Now, for $\e>0$, there exists a $\xi(\cd)\in\BBXi_R[t,T]$ (with large $R>0$) such that %
$$\e+V(\t_0,t,x)\ges J(\t_0,t,x;\xi(\cd))=\dbE\[h(X(T))+\sum_{j\ges1}\ell(\t_j,\xi_j)
+\int_t^T\sum_{j\ges0}g(\t_j,\th,X(\th))d\th\].$$
For this $\xi(\cd)\in\BBXi_R[t,T]$ of form \rf{xi}, define $\xi'(\cd)\in\BBXi_R[t',T]$ as follows
\bel{xi'}\xi'(s)=\sum_{\t_j\les t'}\xi_j{\bf1}_{[t',T]}(s)+\sum_{\t_j>t'}\xi_j{\bf1}_{[\t_j,T]}(s),\qq s\in[t',T].\ee
Thus, this $\xi'(\cd)$ is obtained from $\xi(\cd)$, by moving all the impulses before $t'$ to the moment $t'$. Now, denote $X'(\cd)=X(\cd\,;\t_0,t',x,\xi'(\cd))$. Then
\bel{V<V}\ba{ll}
\ns\ds\e+V(\t_0,t,x)\ges J(\t_0,t,x;\xi(\cd)\\
\ns\ds=\dbE\[h(X(T))\1n+\1n\sum_{\t_j\les t'}\ell(\t_j,\xi_j)\1n+\1n\sum_{\t_j>t'}\ell(\t_j,\xi_j)\1n+\2n\int_t^{t'}\3n
\sum_{\t_j\les
t'}g(\t_j,\th,X(\th))ds\1n+\2n\int_{t'}^T\3n\sum_{\t_j>t'}g(\t_j,\th,X(\th))d\th\]\\
\ns\ds\ges\dbE\[h(X'(T))+\sum_{\t_j\les t'}\ell(t',\xi_j)+\sum_{\t_j>t'}\ell(\t_j,\xi_j)+\int_{t'}^T\sum_{\t_j>t'}g(\t_j,\th,X'(\th))d\th\\
\ns\ds\qq+h(X(T))\1n-\1n h(X'(T))\1n+\2n\int_t^{t'}\3n\sum_{\t_j\les t'}g(\t_j,\th,X(\th))ds\1n+\2n\int_{t'}^T\3n\sum_{\t_j>t'}\1n\(g(\t_j,\th,X(\th))\1n
-\1n g(\t_j,\th,X'(\th))\)d\th\]\\
\ns\ds\ges J(\t_0,t',x;\xi'(\cd))\1n-\1n\dbE\[K|X(T)\1n-\1n X'(T)|\1n+\2n\int_t^{t'}\3n\sum_{\t_j\les t'}|g(\t_j,\th,X(\th))|d\th\1n+\2n\int_{t'}^T\3n\sum_{\t_j>t'}K|X(\th)\1n-\1n X'(\th)|d\th\]\ea\ee
Note that
\bel{state'}\ba{ll}
\ns\ds X'(s)=x+\int_{t'}^sb(\t_0,\th,X'(\th))d\th+\int_{t'}^s\si(\t_0,\th,X'(\th))
dW(\th)\\
\ns\ds\qq\qq+\sum_{\t_j\les t'}\(\int_{t'}^sb(t',\th,X'(\th))d\th
+\int_{t'}^s\si(t',\th,X'(\th))dW(\th)\){\bf1}_{[t',T]}(s)\\
\ns\ds\qq\qq+\sum_{\t_j>t'}\(\int_{\t_j}^sb(\t_j,\th,X'(\th))d\th
+\int_{\t_j}^s\si(\t_j,\th,X'(\th))dW(\th)\){\bf1}_{[\t_j,T]}(s)+\xi'(s),\q s\in[t',T],\ea\ee
Hence, for $s\in[t',T]$,
$$\ba{ll}
\ns\ds\dbE\[\sup_{r\in[t',s]}|X(r)-X'(r)|^p\]=\sum_{j\ges0}\dbE\(\int_{t'}^s
|b(\t_j,\th,X(\th))-b(\t_j,\th,X'(\th))|d\th|\)^p\\
\ns\ds\qq\qq\qq\qq\qq\qq\qq+\sum_{j\ges0}\dbE\[\sup_{r\in[t',s]}\Big|\int_{t'}^r\(\si(\t_j,\th,X(\th))
-\si(\t_j,\th,X'(\th))\)dW(\th)\Big|^p\]\\
\ns\ds\les K\dbE\[\(\int_{t'}^s|X(\th)-X'(\th)|d\th\)^p+\(\int_{t'}^s|X(\th)-X'(\th)|^2d\th\)^{p\over2}
\].\ea$$
This leads to
\bel{X-X'}\dbE\[\sup_{r\in[t',T]}|X(r)-X'(r)|^p\]\les K\dbE|X(t')-X'(t')|^p.\ee
Note
$$\ba{ll}
\ns\ds\dbE\[|X(t')-X'(t')|^p\]=\dbE\Big|\sum_{j\ges0}\(\int_t^{t'}b(\t_j,\th,X(\th))d\th
+\int_t^{t'}\si(\t_j,\th,X(\th)))dW(\th)\)\Big|^p\les K|t-t'|^{p\over2}.\ea$$
Also,
$$\dbE\int_t^{t'}|g(\t_j,\th,X(\th))|d\th\les K|t-t'|.$$
Hence, we obtain
\bel{V>V}\e+V(\t_0,t,x)\ges J(\t_0,t',x)-K|t-t'|^{1\over2}\ges V(\t_0,t',x)-K|t-t'|^{1\over2}.\ee
Now, for $\e>0$, there exists a $\xi'(\cd)\in\BBXi_R[t',T]$ (with large $R$) such that
$$\e+V(\t_0,t',x)\ges J(\t_0,t',x;\xi'(\cd))=\dbE\[h(X'(T))+\sum_{j\ges1}\ell(\t_j,\xi'_j)
+\int_{t'}^T\sum_{j\ges0}g(\t_j,\th,X'(\th))d\th\],$$
where
$$\xi'(s)=\sum_{j\ges1}\xi_j'{\bf1}_{[\t_j',T]}(s),\qq s\in[t',T].$$
For this $\xi'(\cd)$, we define
\bel{xi*}\xi(s)=\sum_{j\ges1}\xi'_j{\bf1}_{[\t_j',T]}(s),\qq s\in[t,T],\ee
Thus, $\xi(\cd)\in\BBXi_R[t,T]$ is an extension of $\xi'(\cd)\in\BBXi_R[t',T]$, from $[t',T]$ to $[t,T]$ by a no-impulse control on $[t,t')$ (a trivial extension). Thus, denoting $X'(\cd)=X(\cd\,;t',x,\xi'(\cd))$, we have
\bel{V<V*}\ba{ll}
\ns\ds\e+V(\t_0,t',x)\ges J(\t_0,x;\xi'(\cd))\\
\ns\ds\ges\dbE\[h(X(T))+\sum_{j\ges1}\ell(\t_j,\xi_j')+\int_t^T\sum_{j\ges0}
g(\t_j,\th,X(\th))d\th\\
\ns\ds\qq+h(X'(T))\1n-\1n h(X(T))\1n-\2n\int_t^{t'}\3n\sum_{j\ges0}g(\t_j,\th,X(\th))ds
\1n+\2n\int_{t'}^T\3n\sum_{\t_j>t'}\2n\(g(\t_j,\th,X'(\th))\1n-\1n g(\t_j,\th,X(\th))\)d\th\]\\
\ns\ds\ges J(\t_0,t,x;\xi(\cd))\1n-\1n\dbE\[K|X'(T)\1n-\1n X(T)|\1n+\2n\int_t^{t'}\3n\sum_{j\ges0}
|g(\t_j,\th,X(\th))|d\th\1n+\2n\int_{t'}^T\3n\sum_{\t_j>t'}\1n K|X'(\th)\1n-\1n X(\th)|d\th\]\\
\ns\ds\ges J(\t_0,t,x;\xi(\cd))-K|t-t'|^{1\over2}\ges V(\t_0,t,x)-K|t-t'|^{1\over2}.\ea\ee
By \rf{V<V} and \rf{V<V*}, we have \rf{V-V} with $\t_0=\t_0'$ and $x=x'$. Hence, \rf{V-V} holds. \endpf

\ms

In what follows, we let $\sV$ be the following:
$$\sV=\Big\{V:\D_*[0,T]\times\dbR^n\to\dbR_+\bigm|V(\cd\,,\cd\,,\cd)\hb{ is continuous and \rf{0<V}--\rf{V-V} hold}\Big\}.$$
Thus, under (H1)--(H2), $V(\cd\,,\cd\,,\cd)\in\sV$. Now, we recall the following:

\bde{semiconvex} Let $G$ be convex. A continuous function $\f:G\to\dbR$ is {\it
semi-convex} if $x\mapsto f(x)-K|x|^2$ is convex for some $K>0$. That is
$$\l\f(x)+(1-\l)\f(x')-\f(\l x+(1-\l)x')\les K\l(1-\l)|x-x'|^2,\qq x,x'\in G,~\l\in[0,1].$$

\ede

Note that linear growth of a function does not necessarily imply it is semi-convex, for example, $\f(x)=\sin(e^x)$, which is not semi-convex. To get the semi-convexity of the value function $V(\cd\,,\cd\,,\cd)$, we need the following additional assumption.

\ms

{\bf(H3)} The map $x\mapsto(b(\t,t,x),\si(\t,t,x))$ is differentiable with
$$|b_x(\t,t,x)-b_x(\t,t,x')|+|\si_x(\t,t,x)-\si_x(\t,t,x')|\les K|x-x'|,\q(\t,t)\in\D_*[0,T],x,x'\in\dbR^n,$$
and the map $x\mapsto(h(x),g(\t,t,x))$ is semi-convex.
\ms

The following gives a further result for the value function.

\bp{semiconvex*} \sl Let {\rm(H1)--(H3)} hold. Then $x\mapsto V(\t_0,t,x)$ is semi-convex, uniformly in $(\t_0,t)\in\D_*[0.T]$.

\ep

\it Proof. \rm Let $(\t_0,t)\in\D_*[0,T]$, and $x,x'\in\dbR^n$. For $\l\in[0,1]$, define
$$x_\l=\l x+(1-\l)x'.\q(\hb{Thus, $x_1=x$, $x_0=x'$})$$
For $\e>0$, there exists a $\xi^\e(\cd)\in\BBXi_R[0,T]$ (with $R>0$ large) such that
$$J(\t_0,t,x_\l;\xi^\e(\cd))<V(\t_0,t,x_\l)+\e.$$
Denote
$$\ba{ll}
\ns\ds X_\l(\cd)=X(\cd\,;\t_0,t,x_\l,\xi^\e(\cd)),\qq\l\in[0,1],\\
\ns\ds X^\l(\cd)=\l X_1(\cd)+(1-\l)X_0(\cd),\qq\l\in[0,1],\ea$$
with
$$X_1(\cd)=X(\cd\,;\t_0,t,x,\xi^\e(\cd))\equiv X(\cd),\qq X_0(\cd)=X(\cd\,;\t_0,t,x',\xi^\e(\cd))\equiv X'(\cd).$$
Then,
$$\ba{ll}
\ns\ds\l V(\t_0,t,x)+(1-\l)V(\t_0,t,x')-V(\t_0,t,x_\l)-\e\\
\ns\ds\les\l J(\t_0,t,x;\xi^\e(\cd))+(1-\l)J(\t_0,t,x';\xi^\e(\cd))-J(\t_0,t,x_\l;\xi^\e(\cd))\\
\ns\ds=\dbE\[\l h(X(T))+(1-\l)h(X'(T))-h(X_\l(T))\\
\ns\ds\qq+\int_t^T\sum_{j\ges0}\(\l g(\t_j,\th,X(\th))+(1-\l)g(\t_j,\th,X'(\th))-g(\t_j,\th,X_\l(\th))\)d\th\]\\
\ns\ds=\dbE\[\l h(X(T))+(1-\l)h(X'(T))-h(X^\l(T))\\
\ns\ds\qq+\int_t^T\sum_{j\ges0}\(\l g(\t_j,\th,X(\th))+(1-\l)g(\t_j,\th,X'(\th))-g(\t_j,\th,X^\l(\th))\)d\th\]\\
\ns\ds\qq+h(X^\l(T))-h(X_\l(T))+\int_t^T\sum_{j\ges0}\(g(\t_j,\th,X_\l(\th))-g(\t_j,\th,X^\l(\th))\)d\th\]\\
\ns\ds\les K\l(1-\l)\dbE\(|X(T)-X'(T)|^2+\int_t^T\sum_{j\ges0}|X(\th)-X'(\th)|^2d\th\)\\
\ns\ds\qq+K\dbE\(|X^\l(T)-X_\l(T)|+\int_t^T\sum_{j\ges0}|X^\l(\th)-X_\l(\th)|d\th\).\ea$$
Since
$$\ba{ll}
\ns\ds\big|\l b(\t,t,x)+(1-\l)b(\t,t,x')-b(\t,t,x_\l)\big|\1n=\1n\big|\l[ b(\t,t,x)-b(\t,t,x_\l)]+(1-\l)[b(\t,t,x')-b(\t,t,x_\l)]\big|\\
\ns\ds=\Big|\l\int_0^1b_x(\t,t,x_\l+\th(1-\l)(x-x'))d\th(1-\l)(x-x')\\
\ns\ds\qq\qq+(1-\l)\int_0^1b_x(\t,t,x_\l+\th\l(x'-x))d\th\l(x'-x)\Big|\\
\ns\ds=\l(1-\l)|x-x'|\,\Big|\2n\int_0^1\3n b_x(\t,t,x_\l+\th(1-\l)(x-x')
-b_x(\t,t,x_\l+\th\l(x'-x))d\th\Big|\les K\l(1-\l)|x-x'|^2.\ea$$
Similarly, we have
$$\big|\l\si(\t,t,x)+(1-\l)\si(\t,t,x')-\si(\t,t,x_\l)\big|\les
K\l(1-\l)|x-x'|^2.$$
Then,
$$\ba{ll}
\ns\ds X^\l(s)=x_\l+\sum_{j\ges0}\(\int_t^s\l b(\t_j,\th,X(\th)){\bf1}_{[\t_j,T]}(\th)d\th
+\int_t^s\l\si(\t_j,\th,X(\th)){\bf1}_{[\t_j,T]}(\th)dW(\th)\)\\
\ns\ds\qq\qq+\sum_{j\ges0}\(\int_t^s(1-\l)b(\t_j,\th,X'(\th)){\bf1}_{[\t_j,T]}(\th)
d\th+\2n\int_t^s(1-\l)\si(\t_j,\th,X'(\th)){\bf1}_{[\t_j,T]}(\th)dW(\th)\)+\xi(s),\\
\ns\ds\qq\qq\qq\qq\qq\qq\qq\qq\qq\qq\qq\qq s\in[t,T],\ea$$
and
$$\ba{ll}
\ns\ds X_\l(s)=x_\l\1n+\2n\sum_{j\ges0}\(\int_t^s\2n b(\t_j,\th,X_\l(\th)){\bf1}_{[\t_j,T]}(\th)
d\th+\2n\int_t^s\2n\si(\t_j,\th,X_\l(\th)){\bf1}_{[\t_j,T]}(\th)dW(\th)\)+\xi(s),\q s\in[t,T],\ea$$
Thus, for any $s\in[t,T]$, by the state equation \rf{state} and Burkholder-Davis-Gundy inequalities,
$$\ba{ll}
\ns\ds\dbE\[\sup_{r\in[t,s]}\big|X^\l(r)-X_\l(r)\big|^2\]\\
\ns\ds+K\sum_{j\ges0}\dbE\[\sup_{r\in[t,s]}\Big|\int_t^r\big[\l b(\t_j,\th,X(\th))
+(1-\l)b(\t_j,\th,X'(\th))-b(\t_j,t,X_\l(\th))\big]{\bf1}_{[\t_j,T]}(\th)d\th\Big|^2\]\\
\ns\ds+K\sum_{j\ges0}\dbE\[\sup_{r\in[t,s]}\Big|\int_t^r\big[\l\si(\t_j,\th,X(\th))
+(1-\l)\si(\t_j,t,X'(\th))-\si(\t_j,\th,X_\l(\th){\bf1}_{[\t_j,T]}(\th))\big]dW(\th)\Big|^2\]\\
\ns\ds+K\sum_{j\ges0}\dbE\[\sup_{r\in[t,s]}\Big|\int_t^r\big[b(\t_j,\th,X^\l(\th))-
b(\t_j,\th,X_\l(\th))\big]{\bf1}_{[\t_j,T]}(\th)d\th\Big|^2\]\\
\ns\ds+K\sum_{j\ges0}\dbE\[\sup_{r\in[t,s]}\Big|\int_t^r\big[\si(\t_j,\th,X^\l(\th))-
\si(\t_j,t,X_\l(\th)){\bf1}_{[\t_j,T]}(\th)\big]dW(\th)\Big|^2\]\\
\ns\ds\les K\l^2(1-\l)^2\dbE\int_t^s|X(\th)-X'(\th)|^4d\th+K\l^2(1-\l)^2\dbE\int_t^s\big|X(\th)
-X'(\th)\big|^4d\th\]\\
\ns\ds\q+K\dbE\[\int_t^s|X^\l(\th)-X_\l(\th)|^2d\th+K\dbE\int_t^s|X^\l(\th)
-X_\l(\th)|^2d\th\]\\
\ns\ds\les K\l^2(1-\l)^2|x-x'|^4+K\int_t^s\dbE|X^\l(\th)-X_\l(\th)|^2d\th.\ea$$
Then, by Gronwall's inequality, we have
$$\dbE\[\sup_{r\in[t,T]}\big|X^\l(r)-X_\l(r)|\]\les K\l(1-\l)|x-x'|^2.$$
Hence,
$$\ba{ll}
\ns\ds\l V(\t_0,t,x)+(1-\l)V(\t_0,t,x')-V(\t_0,t,x_\l)-\e\\
\ns\ds\les\l J(\t_0,t,x;\xi^\e(\cd))+(1-\l)J(\t_0,t,x';\xi^\e(\cd))-J(\t_0,t,x_\l;\xi^\e(\cd))
\les K\l(1-\l)|x-x'|^2.\ea$$
Sending $\e\to0$, we see that $x\mapsto V(\t_0,t,x)$ is semi-convex. \endpf

\section{Dynamic Programming Principle and HJB Equation}

In this section, we are going to derive Bellman's dynamic programming principle (DPP, for short) and the HJB equation, which is a quasi-variational inequality.

\bp{} \sl Let {\rm(H1)--(H2)} hold. Then for any $(\t,t,x)\in\D_*[0,T]\times\dbR^n$,
\bel{V<*}V(\t,t,x)\les\min\Big\{\dbE\[V(\t,t',X(t'))+\int_t^{t'}g(\t,\th,X(\th))
d\th\],N[V(\t,t,\cd)](x)\Big\},\ee
where $X(\cd)=X(\cd\,;\t,t,x,\xi(\cd))$ with $\xi(\cd)\in\BBXi[t,T]$ no impulses in $[t,t']$ for any $t'\in(t,T]$, and
\bel{N}N[V(\t,t,\cd)](x)=\inf_{\xi\in\sK}\[V(\t,t,x)+\ell(t,x+\xi)\].\ee
The operator $N:\sV\to\sV$, and it is Lipschitz continuous:
\bel{N-N}|N[V(\t,t,\cd)](x)-N[\h V(\t,t,\cd)](x)|\les|V(\t,t,x)-\h V(\t,t,x)|,\qq\forall V(\cd\,,\cd,,\cd),\h V(\cd\,,\cd,,\cd)\in\sV.\ee\
If it holds
\bel{V<**}V(\t,t,x)<N[V(\t,t,\cd)](x).\ee
then there exists $\d>0$ such that for all $t'\in[t,t+\d]$,
\bel{V=V}V(\t,t,x)=V(\t,t+\d,X(t+\d))+\int_t^{t+\d}g(\t,s,X(s))ds.\ee
Hence, in the case that $V(\cd\,,\cd\,,\cd)$ is smooth in $(t,x)$, we have
\bel{HJB}\left\{\2n\ba{ll}
\ns\ds\min\Big\{V_t(\t,t,x)+H(\t,t,V_x(\t,t,x),V_{xx}(\t,t,x)),N[V(\t,t,\cd)](x)
-V(\t,t,x)\Big\}=0,\\
\ns\ds\qq\qq\qq\qq\qq\qq\qq\qq\qq\qq\qq\qq(t,x)\in[\t,T]\times\dbR^n,\\
\ns\ds V(\t,T,x)=\min\big\{h(x),N[h](x)\big\},\qq x\in\dbR^n.\ea\right.\ee
where
$$H(\t,t,x,p^\top\2n,P)=pb(\t,t,x)+{1\over2}\si(\t,t,x)^\top P\si(\t,t,x)+g(\t,t,x),\q(t,x,p,P)\in[\t,T]\times\dbR^n\times\dbR^n\times\dbS^n.$$

\ep

\it Proof. \rm Let $(\t,t)\in\D_*[0,T]$ be given. Then for any $t'\in[t,T]$, let $\xi(\cd)\in\BBXi[t,T]$ be an impulse control with no impulses in $[t,t']$, we have
$$V(\t,t,x)\les J(\t,t,x;\xi(\cd))=\dbE\[J(\t,t',X(t');\xi(\cd)\big|_{(t',T]})+\int_t^{t'}g(\t,\th,
X(\th))d\th\].$$
Then, we have
$$V(\t,t,x)\les\dbE\[V(\t,t',X(t'))+\int_t^{t'}g(\t,\th,X(\th))d\th\].$$
On the other hand, if at $t$, we make an impulse in the control, then
$$V(\t,t,x)\les V(\t,t,x+\xi)+\ell(t,\xi).$$
Thus, we must have
$$V(\t,t,x)\les\inf_{\xi\in\sK}\(V(\t,t,x+\xi)+\ell(t,\xi)\)\equiv N[V(\t,t,\cd)](x).$$
Combining the above, we have \rf{V<*}. Further, by the continuity of $V(\cd\,,\cd\,,\cd)$ and $\ell(\cd\,,\cd)$, we have
$$|V(\t,t,x+\xi)+\ell(t,\xi)-[V(\t',t',x'+\xi)+\ell(t',\xi)]|\les\o(|\t-\t'|+|t-t'|)
+K\big(|t-t'|^{1\over2}+|x-x'|\big).$$
Thus, by the definition of operator $N$, we have
$$|N[V(\t,t,\cd)](x)-N[V(\t',t',\cd)](x')|\les\o\big(|\t-\t'|+|t-t'|\big)+K\big(
|t'-t|^{1\over2}+|x-x'|\big).$$
Hence, $N:\sV\to\sV$ and we see obviously that \rf{N-N}.

\ms

Next, if \rf{V<**} holds, by the continuity of $N[V(\t,t,\cd)](x)$, we see that
on some interval $[t,t+\d]$, any impulse control containing impulses in this interval will not be optimal. Hence,
$$\ba{ll}
\ns\ds V(\t,t,x)=J(\t,t,x;\xi(\cd))=\dbE\[J(\t,t',X(t');\xi(\cd))+\int_t^{t'}g(\t,\th,
X(\th))d\th\]\\
\ns\ds\qq\qq\qq\qq\qq\q\ges\dbE\[V(\t,t',X(t'))+\int_t^{t'}g(\t,\th,
X(\th))d\th\],\ea$$
with $\xi(\cd)$ contains no impulses in $[t,t']$. Then there exists $\d>0$ such that \rf{V=V} holds. Hence, if $(t,x)\mapsto V(\t,t,x)$ is smooth, then by It\^o's formula, we see that $V(\cd\,,\cd\,,\cd)$ satisfies the equation in \rf{HJB}. There are two cases: In the case that
$$h(x)<N[h](x)\equiv\inf_{\xi\in\sK}\(h(x+\xi)+\ell(T,\xi)\).$$
Then, at $T$, optimal impulse control will have no impulses. Hence,
$$V(\t,T,x)=h(x),\qq x\in\dbR^n.$$
In the other case, i.e.,
$$h(x)=N[h](x).$$
Then, at $T$, there would be an impulse so that
$$N[h](x)=h(x+\xi)+\ell(T,\xi).$$
Then we must have
$$V(\t,T,x)=N[h](x).$$
Hence, \rf{HJB} holds. \endpf

\ms

Next result is interesting.

\bp{3.2} \sl Let {\rm(H1)--(H2)} hold. Let $(\t,t,x)\in\D_*[0,T]\times\dbR^n$ so that
\bel{V=N}V(\t,t,x)=N[V(\t,x,\cd)](x)\equiv\inf_{\xi\in\sK}\(V(\t,t,x+\xi)
+\ell(t,\xi)\)=V(\t,t,x+\h\xi)+\ell(t,\h\xi),\ee
for some $\h\xi\in\sK$. Then
\bel{V<N}V(\t,t,x+\h\xi)<N[V(\t,t,\cd)](x+\h\xi).\ee

\ep

\it Proof. \rm Let $\h\xi\in\sK$ such that \rf{V=N} holds. Then
we claim that \rf{V<N} holds.

In fact, let
$$V(\t,t,x+\h\xi\,)=N[V(\t,t,\cd)](x+\h\xi)=V(\t,t,x+\h\xi+\wt\xi\,)
+\ell(t,\wt\xi\,),$$
for some $\wt\xi\in\sK$. Then
$$\ba{ll}
\ns\ds V(\t,t,x)=V(\t,t,x+\h\xi)+\ell(t,\h\xi)=V(\t,t,x+\h\xi+\wt\xi\,)+\ell(t,\wt\xi\,)
+\ell(t,\h\xi)\\
\ns\ds\qq\qq\qq\qq\qq\qq\qq\q~>V(\t,t,x+\h\xi+\wt\xi\,)+\ell(t,\h\xi+\wt\xi\,)\ges
N[V(\t,t,\cd)](x),\ea$$
which contradicts \rf{ell}. Thus, \rf{V<N} has to be true. \endpf

\ms

The above result implies if
\bel{bar xi}\bar\xi(\cd)=\sum_{j\ges1}\bar\xi_j{\bf1}_{[\bar\t_j,T]}(\cd)\in
\BBXi[t,T]\ee
is an optimal impulse control, then
\bel{t<t}t\les\bar\t_1<\bar\t_2<\cds<\bar\t_{\k(\bar\xi(\cd))}\les T.\ee
In particular, from now on, for any $\xi(\cd)\in\BBXi[t,T]$ of form \rf{xi}, we always assume that
\bel{t<t*}t\les\t_1<\t_2<\cds<\t_{\k(\xi(\cd))}\les T.\ee

\ms

We have seen that \rf{HJB} is an HJB equation on $[\t,T]\times\dbR^N$ with $\t\in[0,T)$ as a parameter. Hence, we have the following natural definition.

\bde{viscosity} (i) Continuous function $V(\cd\,,\cd\,,\cd)$ is {\it viscosity
subsolution} (resp. viscosity supersolution) of \rf{HJB} if
$$V(\t,T,x)\les\min\{h(x),N[h](x)\},\q(V(\t,T,x)\ges\min\{h(x),N[h](x))$$
and if $V(\t,t,x)-\f(\t,t,x)$ attains a local maximum (resp. local minimum) at
$(t_0,x_0)$, then
$$\ba{ll}
\ns\ds\min\Big\{\f_t(\t,t_0,x_0)+H(\t,t_0,x_0,\f_x(\t,t_0,x_0),\f_{xx}(\t,
t_0,x_0)),N[V(\t,\cd\,,\cd)](t_0,x_0)-V(\t,t_0,x_0)\Big\}\ges0\\
\ns\ds\qq\qq\qq\qq\qq\qq\qq\qq\qq\qq\qq\qq\qq\qq\qq\qq\qq\qq\qq\qq(\les0)\ea$$

(ii) Continuous function $V(\cd\,,\cd\,,\cd)$ is a viscosity solution of \rf{HJB}
if it is both viscosity subsolution and viscosity supersolution of \rf{HJB}.
\ede

Now, we state the following result whose proof follows some standard arguments (see
\cite{Crandall-Lions 1983, Crandall-Ishii-Lions 1992, Tang-Yong 1993, Yong 1994, Li-Yong 2021}).

\bp{} \sl Let {\rm(H1)--(H2)} hold. Then the value function $V(\t,t,x)$ is the unique viscosity solution of \rf{HJB}.

\ep

\section{Optimal Impulse Controls and Maximum Principle}

We have seen that any optimal impulse control of our Problem (IC) can only have finitely many impulses. Thus, we have the following.

\bt{existence} \sl Let {\rm(H1)--(H2)} hold. Then, for given $(\t,t,x)\in\D_*[0,T]\times\dbR^n$, Problem {\rm(IC)} admits an optimal impulse control.

\et

\it Proof. \rm Let $(\t,t,x)\in\D_*[0,T]\times\dbR^n$. Let $\xi^k(\cd)\in\BBXi[t,T]$ be a minimizing sequence, i.e.,
$$J(\t,t,x;\xi^k(\cd))\to V(\t,t,x).$$
Clearly, we may assume that each $\xi^k(\cd)$ has the same number of impulses. Thus,
$$\xi^k(s)=\sum_{j\ges1}\xi^k_j{\bf1}_{[\t_j^k,T]}(s),\qq s\in[t,T],~k\ges1.$$
Then, by extracting sub-sequences, we see the existence of an optimal control.\endpf

\ms

To obtain a necessary condition of Pontryagin type, we will have three standard steps: Derive Taylor type expansion, obtain a variation inequality; and use duality to get maximum condition. We now make these more precise.

\subsection{Taylor Expansion}

We first make the following additional assumption.

\ms

{\bf(H4)} The map $x\mapsto(b(\t,s,x),\si(\t,s,x),g(\t,s,x),h(x))$ is twice continuously differentiable, and the map $\t\mapsto(b(\t,s,x),\si(\t,s,x),\ell(\t,\xi))$ is continuously differentiable.

\ms

The following is classical.

\bl{Taylor} \sl Let $\f:\dbR^n\to\dbR$ be continuous twice differentiable. Then for any $x_0\in\dbR^n$,
\bel{f=}\f(x)=\f(x_0)+\f_x(x_0)(x-x_0)+(x-x_0)^\top\(\int_0^1\int_0^1\f_{xx}
(x_0+\a\b(x-x_0)))\b d\b d\a\)(x-x_0)\ee
for any $x\in\dbR^n$.
\el

\it Proof. \rm By Cauchy type Taylor expansion, we have
$$\ba{ll}
\ns\ds\f(x)=\f(x_0)+\(\int_0^1\f_x(x_0+\a(x-x_0))d\a\)(x-x_0)\\
\ns\ds\qq=\f(x_0)+\f_x(x_0)(x-x_0)+\(\int_0^1\big[\f_x(x_0+\a(x-x_0)))
-\f_x(x_0)\big]d\a\)(x-x_0)\\
\ns\ds\qq=\f(x_0)+\f_x(x_0)(x-x_0)+(x-x_0)^\top\(\int_0^1\int_0^1\f_{xx}
(x_0+\a\b(x-x_0)))\b d\b d\a\)(x-x_0)\ea$$
Then \rf{f=} follows. \endpf

\ms

In the case that $\f:\dbR^n\to\dbR^m$ is smooth, we have a similar result. In this case, we denote
$$\ba{ll}
\ns\ds(x-x_0)^\top\(\int_0^1\3n\int_0^1\f_{xx}
(x_0+\a\b(x-x_0)))\b d\b d\a\)(x-x_0)\\
\ns\ds=(x-x_0)^\top\begin{pmatrix}\ds\int_0^1\3n\int_0^1\f^1_{xx}
(x_0+\a\b(x-x_0)))\b d\b d\a\\
\ds\int_0^1\3n\int_0^1\f^2_{xx}
(x_0+\a\b(x-x_0)))\b d\b d\a\\
\vdots\\
\ds\int_0^1\3n\int_0^1\f^m_{xx}
(x_0+\a\b(x-x_0)))\b d\b d\a\end{pmatrix}(x-x_0),\ea$$
which is an $\dbR^m$ vector.

\ms

Next result gives Taylor expansion about the optimal pair $(\bar X(\cd),\bar\xi(\cd))$.

\bt{Taylor Expansion} \sl Let {\rm(H1)--(H2)} and {\rm(H4)} hold. Let $(\t_0,t,x)\in\D_*[0,T]\times\dbR^n$ be given and $\bar\xi(\cd)\in\BBXi[t,T]$ be an optimal impulse control of the following form \rf{bar xi}
with the optimal state process being $\bar X(\cd)=X(\cd\,;\t_0,t,x,\xi(\cd))$. Let $i\in\{1,2,\cds,\k(\bar\xi(\cd))\}$, $\e,\bar\e\in(0,1)$ small, $\eta\in\sK$, and
\bel{xi^e}\ba{ll}
\ns\ds\xi^{i,\e,\bar\e}(\cd)=\sum_{j\ges1, j\ne i}\bar\xi_j{\bf1}_{[\bar\t_j,T]}(\cd)+[\bar\xi_i+\e(\eta-\bar\xi_i)]
{\bf1}_{[\bar\t_i+\bar\e,T]}(\cd)\\
\ns\ds\qq\q~=\bar\xi(\cd)-\bar\xi_i{\bf1}_{[\bar\t_i,
\bar\t_i+\bar\e)}(\cd)+\e(\eta-\bar\xi_i){\bf1}_{[\bar\t_i+\bar\e,T]}\in\BBXi[t,T],
\ea\ee
be the $i$-th perturbation of $\bar\xi(\cd)$. Let $X^{i,\e,\bar\e}(\cd)=X(\cd\,;\t_0,t,x,\xi^{i,\e,\bar\e}(\cd))$ be the corresponding process. Let
\bel{X_1}\ba{ll}
\ns\ds X^{i,\e,\bar\e}_1(s)=\int_t^sB_x(\th)X^{i,\e,\bar\e}_1(\th)d\th+\int_t^s\(\Si_x(\th)
X^{i,\e,\bar\e}_1(\th)-\si(\bar\t_i,\th){\bf1}_{[\bar\t_i,
\bar\t_i+\bar\e)}(\th)\)dW(\th)\\
\ns\ds\qq\qq\qq\qq-\bar\xi_i{\bf1}_{[\bar\t_i,\bar\t_i+\bar\e)}(s),\qq s\in[t,T],\ea\ee
and
\bel{X_2}\ba{ll}
\ns\ds X^{i,\e,\bar\e}_2(s)=\int_t^s\(B_x(\th)X^{i,\e,\bar\e}_2(\th)+
X^{i,\e,\bar\e}_1(\th)^\top B_{xx}(\th)X_1^{i,\e,\bar\e}(\th)-b(\bar\t_i,\th){\bf1}_{[\bar\t_i,\bar\t_i+\bar\e)}
(\th)\)d\th\\
\ns\ds\qq\qq\qq+\int_t^s\(\Si_x(\th)X^{i,\e,\bar\e}_2(\th)+X^{i,\e,\bar\e}_1
(\th)^\top\Si_{xx}(\th)
X^{i,\e,\bar\e}_1(\th)\)dW(\th)\\
\ns\ds\qq\qq\qq+\bar\e\z^i(s)+\e\xi_i{\bf1}_{[\bar\t_i,T]}(s),\qq s\in[t,T],\ea\ee
where
\bel{bB}\left\{\2n\ba{ll}
\ds b(\bar\t_i,\th)=b(\bar\t_i,\th,\bar X(\th)),\qq \si(\bar\t_i,\th)=\si(\bar\t_i,\th,\bar X(\th)),\q i\ges0,\\ [2mm]
\ns\ds b^i_\t(\th)=b_\t(\bar\t_i,\th,\bar X(\th)){\bf1}_{[\bar\t_i,T]}(\th),\q\si^i_\t(\th)=\si_\t(\bar\t_i,\th,\bar X(\th)){\bf1}_{[\bar\t_i,T]}(\th),\q i\ges0,\\ [2mm]
\ns\ds\z^i(s)=\int_t^sb_\t^i(\th)d\th+\int_t^s\si_\t^i(\th)dW(\th)\\
\ns\ds\qq\equiv\int_t^sb_\t(\bar\t_i,\th,\bar X(\th){\bf1}_{[\bar\t_i,T]}(\th)d\th+\int_t^s\si_\t(\bar\t_i,\th,\bar X(\th)){\bf1}_{[\bar\t_i,T]}(r)dW(\th),\\
\ns\ds B_x(\th)=\sum_{j\ges0}b_x(\bar\t_j,\th,\bar X(\th)){\bf1}_{[\bar\t_j,T]}(\th),\q\Si_x(\th)=\sum_{j\ges0}\si_x(\bar\t_j,\th,\bar X(\th)){\bf1}_{[\bar\t_j,T]}(\th),\\
\ns\ds B_{xx}(\th)={1\over2}\sum_{j\ges0}b_{xx}(\bar\t_j,\th,\bar X(\th)){\bf1}_{[\bar\t_j,T]}(\th),\q\Si_{xx}(\th)={1\over2}\sum_{j\ges0}\si_{xx}(\bar\t_j,\th,\bar X(\th)){\bf1}_{[\bar\t_j,T]}(\th),\ea\right.\ee
with $\bar\t_0=\t_0$. Then for $m\ges1$,
\bel{|X-X|,|X_1|}\dbE\[\int_t^T|X^{i,\e,\bar\e}(r)-\bar X(r)|^{2m}dr\]=O(\e^{2m}+\bar\e),\q\dbE\[\int_t^T|X_1^{i,\e,\bar\e}(r)|^{2m}dr\]=O(\bar\e),\ee
\bel{|X-X-X_1|,|X_2|}\dbE\[\int_t^T|X^{i,\e,\bar\e}(r)-\bar X(r)-X^{i,\e,\bar\e}_1(r)|^{2m}dr\]=O(\e^{2m}+\bar\e^{2m}),\q \dbE\[\int_t^T|X^{i,\e,\bar\e}_2(r)|^{2m}dr\]=O(\e^{2m}+\bar\e^{2m}),\ee
\bel{|X-X-X-X|}\dbE\[\int_t^T|X^{i,\e,\bar\e}(r)-\bar X(r)-X^{i,\e,\bar\e}_1(r)-X^{i,\e,\bar\e}_2(r)|^{2m}dr\]
=o(\e^{2m}+\bar\e^{2m}).\ee

\et

Now, let us make some comments on the above result. First of all, by Proposition \ref{3.2}, we see that
$$t\les\bar\t_1<\bar\t_2<\cds<\bar\t_{\k(\bar\xi(\cd))}\les T.$$
Thus, perturbation \rf{xi^e} is well-defined. Next, by the convexity of $\sK$, for any $\eta\in\sK$,
$$\bar\xi_j+\e(\eta-\bar\xi_j)=(1-\e)\bar\xi_j+\e\eta\in\sK.$$
Hence, the inclusion in \rf{xi^e} holds. In what follows, for convenience, we will denote $\xi_i=\eta-\bar\xi_i$. Note that in order that the inclusion in \rf{xi^e} holds, we have assumed $\bar\t_{\k(\bar\xi(\cd))}<T$. But $\bar\t_1=t$ is allowed.

\ms

Note that \rf{|X-X-X-X|} is a Taylor type expansion. In \cite{Hu-Yong 1991}, the variation was made with respect to the sizes $\bar\xi_j$ of the optimal impulse control only (see \cite{Yong-Zhang 1992} also). The variations with respect to the optimal impulse control moments is hard to apply since $X_1^{i,\e,\bar\e}(\cd)$ does not have a limit, in general. On the other hand, the idea of using a Taylor expansion and the first and the second adjoint equations \cite{Peng 1990} (see also \cite{Yong-Zhou 1999, Yong 2022}) seems to be helpful. However, the following might not be small (for $\e,\bar\e\to0$):
$$\dbE\[\sup_{r\in[t,T]}|X^{i,\e,\bar\e}(r)-\bar X(r)-X_1^{i,\e,\bar\e}(r)-X_2^{i,\e,\bar\e}(r)|^{2m}\],$$
and moreover, each term in the variation of the cost functional will be different. Thus, we need to make creative modifications of existing techniques to achieved what we expected.

\ms

Also, in the case $\bar\t_1>t$ (no restrictions on the other $\bar\t_j$, including the case $\bar\t_{\k(\bar\xi(\cd))}=T$), we may replace
$\bar\e$ by $-\bar\e$. Thus, we have
\bel{xi-}\ba{ll}
\ns\ds\xi^{i,\e,-\bar\e}(\cd)=\sum_{j\ges1,j\ne i}\bar\xi_j{\bf1}_{[\bar\t_j,T]}(\cd)+[(1-\e)\bar\xi_i+\e\eta]{\bf1}_{[\bar\t_i
-\bar\e,T]}(\cd)\\
\ns\ds\qq\q~=\sum_{j\ges1}\bar\xi_j{\bf1}_{[\bar\t_j,T]}(\cd)+\e\xi_i{\bf1}_{[\bar\t_i
-\bar\e,T]}(\cd)-\bar\xi_i
{\bf1}_{[\bar\t_i-\bar\e,\bar\t_i)}(\cd),\ea\ee
In what follows, we are going to consider a perturbation of form \rf{xi^e} only
(the derivation of the perturbation \rf{xi-} will be similar).

\ms

\it Proof. \rm For given $(\t_0,t,x)\in\D_*[0,T]\times\dbR^n$, let $\bar\xi(\cd)$ be an optimal impulse control of form \rf{bar xi}. Then (recall $\bar\t_0=\t_0$)
$$\ba{ll}
\ns\ds \bar X(s)=x+\sum_{j\ges0}\(\int_t^sb(\bar\t_j,\th,\bar X(\th)){\bf1}_{[\bar\t_j,T]}(\th)d\th
+\2n\int_t^s\2n\si(\bar\t_j,\th,\bar X(\th)){\bf1}_{[\bar\t_j,T]}(\th)dW(\th)\)\1n+\1n\sum_{j\ges1}\bar\xi_j
{\bf1}_{[\bar\t_j,T]}(s),\\
\ns\ds\qq\qq\qq\qq\qq\qq\qq\qq\qq\qq\qq\qq\qq s\in[t,T],\ea$$
Also,
$$\ba{ll}
\ns\ds X^{i,\e,\bar\e}(s)=X(s;\t_0,t,x,\xi^{i,\e,\bar\e}(\cd))\\
\ns\ds= x+\int_t^s\(\sum_{j\ges0,j\ne i}b(\bar\t_j,\th,X^{i,\e,\bar\e}(\th)){\bf1}_{[\bar\t_j,T]}(\th)
+b(\bar\t_i+\bar\e,\th,X^{i,\e,\bar\e}(\th)){\bf1}_{[\bar\t_i+\bar\e,T]}(\th)\)d\th\\
\ns\ds\qq\qq+\3n\sum_{j\ges1,j\ne i}\bar\xi_j{\bf1}_{[\t_j,T]}(s)+(\bar\xi_i+\e\xi_i){\bf1}_{[\bar\t_i+\bar\e,T]}(s)\\
\ns\ds\qq\qq+\2n\int_t^s\2n\(\2n\sum_{j\ges0,j\ne i}\si(\bar\t_j,\th,X^{i,\e,\bar\e}(\th)){\bf1}_{[\bar\t_j,T]}(\th)+\si(\bar\t_i+\bar\e,
\th,X^{i,\e,\bar\e}(\th)){\bf1}_{[\bar\t_i+\bar\e,T]}(\th)\)dW(\th)\\
\ns\ds=x+\2n\sum_{j\ges0}\(\int_t^s\3n b(\bar\t_j,\th,X^{i,\e,\bar\e}(\th)){\bf1}_{[\bar\t_i,T]}(\th)
d\th+\2n\int_t^s\3n\si(\bar\t_j,\th,X^{i,\e,\bar\e}(\th)){\bf1}_{[\bar\t_i,T]}(\th)d\th\)
+\sum_{j\ges1}\bar\xi_j{\bf1}_{[\bar\t_i,T]}(s)\\
\ns\ds\qq\qq+\int_t^s\(b(\bar\t_i+\bar\e,\th,X^{i,\e,\bar\e}(\th)){\bf1}_{[\bar\t_i
+\bar\e,T]}(\th)-b(\bar\t_i,\th,X^{i,\e,\bar\e}(\th))
{\bf1}_{[\bar\t_i,T]}(\th)\)d\th\\
\ns\ds\qq\qq+\int_t^s\(\si(\bar\t_i+\bar\e,\th,X^{i,\e,\bar\e}(\th)){\bf1}_{[\bar\t_i
+\bar\e,T]}(\th)-\si(\bar\t_i,\th,
X^{i,\e,\bar\e}(\th)){\bf1}_{[\bar\t_i,T]}(\th)\)dW(\th)\\
\ns\ds\qq\qq+\e\xi_i{\bf1}_{[\bar\t_i+\bar\e,T]}(s)+\bar\xi_i
[{\bf1}_{[\bar\t_i+\bar\e,T]}(s)-{\bf1}_{[\bar\t_i,T]}(s)]\\
\ns\ds=x+\2n\sum_{j\ges0}\(\int_t^s\3n b(\bar\t_j,\th,X^{i,\e,\bar\e}(\th)){\bf1}_{[\bar\t_i,T]}(\th)
d\th+\2n\int_t^s\3n\si(\bar\t_j,\th,X^{i,\e,\bar\e}(\th)){\bf1}_{[\bar\t_i,T]}(\th)d\th\)
+\sum_{j\ges1}\bar\xi_j{\bf1}_{[\bar\t_i,T]}(s)\\
\ns\ds\qq\q+\2n\int_t^s\3n\(\int_0^1b_\t(\bar\t_i+\a\bar\e,\th,X^{i,\e,\bar\e}(\th))
d\a\)\bar\e{\bf1}_{[\bar\t_i
+\bar\e,T]}(\th)-b(\bar\t_i,\th,X^{i,\e,\bar\e}(\th))
{\bf1}_{[\bar\t_i,\bar\t_i+\bar\e)}(\th)\)d\th\\
\ns\ds\qq\q+\2n\int_t^s\3n\(\int_0^1\si_\t(\bar\t_i+\a\bar\e,\th,X^{i,\e,\bar\e}(\th))d\a\)
\bar\e{\bf1}_{[\bar\t_i+\bar\e,T]}(\th)-\si(\bar\t_i,\th,
X^{i,\e,\bar\e}(\th)){\bf1}_{[\bar\t_i,\bar\t_i+\bar\e)}(\th)\)dW(\th)\\
\ns\ds\qq\qq+\e\xi_i{\bf1}_{[\bar\t_i+\bar\e,T]}(s)-\bar\xi_i
[{\bf1}_{[\bar\t_i,\bar\t_i+\bar\e)}(s)\\
\ns\ds=x+\2n\sum_{j\ges0}\(\int_t^s\3n b(\bar\t_j,\th,X^{i,\e,\bar\e}(\th)){\bf1}_{[\bar\t_i,T]}(\th)
d\th+\2n\int_t^s\3n\si(\bar\t_j,\th,X^{i,\e,\bar\e}(\th)){\bf1}_{[\bar\t_i,T]}(\th)
dW(\th)\)+\sum_{j\ges1}\bar\xi_j{\bf1}_{[\bar\t_i,T]}(s)\\
\ns\ds\qq\q-\int_t^sb(\bar\t_i,\th,X^{i,\e,\bar\e}(\th))
{\bf1}_{[\bar\t_i,\bar\t_i+\bar\e)}(\th)d\th-\2n\int_t^s\3n\si(\bar\t_i,\th,
X^{i,\e,\bar\e}(\th)){\bf1}_{[\bar\t_i,\bar\t_i+\bar\e)}(\th)dW(\th)\\
\ns\ds\qq\q+\bar\e\(\int_t^sb_\t^{i,\e,\bar\e}(\th)d\th+\int_t^s\3n\si_\t^{i,\e,\bar\e}(\th)
dW(\th)\)+\e\xi_i{\bf1}_{[\bar\t_i+\bar\e,T]}(s)-\bar\xi_i
[{\bf1}_{[\bar\t_i,\bar\t_i+\bar\e)}(s),\ea$$
where
\bel{b to}\ba{ll}\left\{\2n\ba{ll}
\ds b^{i,\e,\bar\e}_\t(\th)=\(\int_0^1b_\t(\bar\t_i+\a\bar\e,\th,X^{i,\e,\bar\e}(\th))
d\a\){\bf1}_{[\bar\t_i+\bar\e,T]}(\th)\to b^i_\t(\th),\\
\ns\ds\si^{i,\e,\bar\e}_\t(\th)=\(\int_0^1\si_\t(\bar\t_i+\a\bar\e,\th,X^{i,\e,\bar\e}(\th))
d\a\){\bf1}_{[\bar\t_i+\bar\e,T]}(\th)\to\si_\t^i(\th),\ea\right.\qq\as\\
\ns\ds\qq\qq\qq\qq\qq\hb{as long as }\lim_{\e,\bar\e\to0}X_1^{i,\e,\bar\e}(\th)=\bar X(\th),\qq\as\ea\ee
Hence,
$$\ba{ll}
\ns\ds X^{i,\e,\bar\e}(s)-\bar X(s)=\sum_{j\ges0}\[\int_t^s\(b(\bar\t_j,\th,X^{i,\e,\bar\e}
(\th))-b(\bar\t_j,\th,\bar X(\th))\){\bf1}_{[\bar\t_i,T]}(\th)d\th\\
\ns\ds\qq\qq\qq\qq\qq+\int_t^s\(\si(\bar\t_j,\th,X^{i,\e,\bar\e}(\th))
-\si(\bar\t_j,\th,\bar X(\th))\){\bf1}_{[\bar\t_i,T]}(\th)dW(\th)\]\\
\ns\ds\qq\qq\qq\qq\qq-\2n\int_t^s\3n b(\bar\t_i,\th,X^{i,\e,\bar\e}(\th))
{\bf1}_{[\bar\t_i,\bar\t_i+\bar\e)}(\th)d\th-\2n\int_t^s\3n\si(\bar\t_i,\th,
X^{i,\e,\bar\e}(\th)){\bf1}_{[\bar\t_i,\bar\t_i+\bar\e)}(\th)dW(\th)\\
\ns\ds\qq\qq\qq\qq\qq+\bar\e\(\int_t^sb_\t^{i,\e,\bar\e}(\th)d\th
+\int_t^s\si_\t^{i,\e,\bar\e}(\th)
dW(\th)\)+\e\xi_i{\bf1}_{[\bar\t_i+\bar\e,T]}(s)-\bar\xi_i
[{\bf1}_{[\bar\t_i,\bar\t_i+\bar\e)}(s).\ea$$
Then, for $m\ges1$,
$$\ba{ll}
\ns\ds\dbE\int_t^s\3n|X^{i,\e,\bar\e}(r)-\bar X(r)|^{2m}dr\les K\Big\{\dbE\int_t^s\[\(\int_t^r\3n|X^{i,\e,\bar\e}(\th)-\bar X(\th)|d\th\)^{2m} dr+\(\int_t^r\3n|X^{i,\e,\bar\e}(\th)-\bar X(\th)|^2d\th\)^m\\
\ns\ds\qq+\(\int_t^r{\bf1}_{[\bar\t_i,\bar\t_i+\bar\e)}(\th)d\th\)^{2m}
+\(\int_t^r{\bf1}_{[\bar\t_i,\bar\t_i+\bar\e)}(\th)d\th\)^m\]dr
+\bar\e^{2m}+\e^{2m}+\int_t^s{\bf1}_{[\bar\t_i,\bar\t_i+\bar\e)}(r)dr\]\\
\ns\ds\les K\[\int_t^s\(\dbE\int_t^r|X^{i,\e,\bar\e}(\th)-\bar X(\th)|^{2m}d\th\)dr+\e^{2m}+\bar\e\].\ea$$
By Gronwall's inequality, we have the first of \rf{|X-X|,|X_1|}, which gives
$$\lim_{\e,\bar\e\to0}X^{i,\e,\bar\e}(s)=\bar X(s),\qq\as$$
On the other hand, by Lemma \ref{Taylor}, we have
$$\ba{ll}
\ns\ds X^{i,\e,\bar\e}(s)-\bar X(s)\\
\ns\ds=\sum_{j\ges0}\[\int_t^s\(b(\bar\t_j,\th,X^{i,\e,\bar\e}(\th))
-b(\bar\t_j,\th,\bar X(\th))\){\bf1}_{[\bar\t_i,T]}(\th)d\th\\
\ns\ds\qq+\int_t^s\(\si(\bar\t_i,\th,X^{i,\e,\bar\e}(\th))-\si(\bar\t_i,\th,\bar X(\th))\){\bf1}_{[\bar\t_i,T]}(\th)dW(\th)\]\\
\ns\ds\qq-\int_t^sb(\bar\t_i,\th,X^{i,\e,\bar\e}(\th))
{\bf1}_{[\bar\t_i,\bar\t_i+\bar\e)}(\th)d\th-\int_t^s\si(\bar\t_i,\th,
X^{i,\e,\bar\e}(\th)){\bf1}_{[\bar\t_i,\bar\t_i+\bar\e)}(\th)dW(\th)\\
\ns\ds\qq+\bar\e\(\int_t^sb_\t^{i,\e,\bar\e}(\th)d\th
+\int_t^s\si_\t^{i,\e,\bar\e}(\th)
dW(\th)\)+\e\xi_i{\bf1}_{[\bar\t_i+\bar\e,T]}(s)-\bar\xi_i
[{\bf1}_{[\bar\t_i,\bar\t_i+\bar\e)}(s)\\
\ns\ds=\int_t^s\sum_{j\ges0}\[b_x(\bar\t_j,\th)[X^{i,\e,\bar\e}(\th)-\bar X(\th)]\\
\ns\ds\q+[X^{i,\e,\bar\e}(\th)\1n-\1n\bar X(\th)]^\top\(\int_0^1\3n\int_0^1\1n b_{xx}(\bar\t_j,\th,\bar X(\th)\1n+\1n\a\b[X^{i,\e,\bar\e}(\th)\1n-\1n\bar X(\th)])\b d\b d\a\)[X^{i,\e,\bar\e}(\th)-\bar X(\th)]\]{\bf1}_{[\bar\t_j,T]}(\th)d\th\\
\ns\ds\q+\int_t^s\sum_{j\ges0}\[\si_x(\bar\t_j,\th)[X^{i,\e,\bar\e}(\th)-\bar X(\th)]\ea$$
$$\ba{ll}
\ns\ds\q+[X^{i,\e,\bar\e}(\th)\1n-\1n\bar X(\th)]^\top\1n\(\2n\int_0^1\3n\int_0^1\3n\si_{xx}(\bar\t_j,\th,
\bar x(\th)\1n+\1n\a\b[X^{i,\e,\bar\e}(\th)\1n-\1n\bar X(\th)])\b d\b d\a\1n\)[X^{i,\e,\bar\e}(\th)\1n-\1n\bar X(\th)]\]{\bf1}_{[\bar\t_j,T]}(\th)dW(\th)\\
\ns\ds\q-
\int_t^sb(\bar\t_i,\th,X^{i,\e,\bar\e}(\th))
{\bf1}_{[\bar\t_i,\bar\t_i+\bar\e)}(\th)d\th-\int_t^s\si(\bar\t_i,\th,
X^{i,\e,\bar\e}(\th)){\bf1}_{[\bar\t_i,\bar\t_i+\bar\e)}(\th)dW(\th)\\
\ns\ds\q+\bar\e\(\int_t^sb_\t^{i,\e,\bar\e}(\th)d\th+\int_t^s\si_\t^{i,\e,\bar\e}(\th)
dW(\th)\)+\e\xi_i{\bf1}_{[\bar\t_i+\bar\e,T]}(s)
-\bar\xi_i{\bf1}_{[\bar\t_i,\bar\t_i+\bar\e)}(s)\ea$$
$$\ba{ll}
\ns\ds=\int_t^s\(B_x(\th)[X^{i,\e,\bar\e}(\th)-\bar X(\th)]+[X^{i,\e,\bar\e}(\th)-\bar X(\th)]^\top B_{xx}^{i,\e,\bar\e}(\th)[X^{i,\e,\bar\e}(\th)-\bar X(\th)]\)d\th\qq\qq\qq\qq\\
\ns\ds\q+\int_t^s\(\Si_x(\th)[X^{i,\e,\bar\e}(\th)-\bar X(\th)]+[X^{i,\e,\bar\e}(\th)-\bar X(\th)]^\top\Si_{xx}^{i,\e,\bar\e}(\th)
[X^{i,\e,\bar\e}(\th)-\bar X(\th)]\)dW(\th)\\
\ns\ds\q-\int_t^sb(\bar\t_i,\th,X^{i,\e,\bar\e}(\th))
{\bf1}_{[\bar\t_i,\bar\t_i+\bar\e)}(\th)d\th-\int_t^s\si(\bar\t_i,\th,
X^{i,\e,\bar\e}(\th)){\bf1}_{[\bar\t_i,\bar\t_i+\bar\e)}(\th)dW(\th)\\
\ns\ds\q+\bar\e\(\int_t^sb_\t^{i,\e,\bar\e}(\th)d\th+\int_t^s\si_\t^{i,\e,\bar\e}(\th)
dW(\th)\)+\e\xi_i{\bf1}_{[\bar\t_i+\bar\e,T]}(s)
-\bar\xi_i{\bf1}_{[\bar\t_i,\bar\t_i+\bar\e)}(s),\ea$$
where
\bel{B to}\ba{ll}
\ns\ds B^{i,\e,\bar\e}_{xx}(\th)=\sum_{j\ges0}\( \int_0^1\3n\int_0^1\1n b_{xx}(\bar\t_j,\th,\bar X(\th)\1n+\1n\a\b[X^{i,\e,\bar\e}(\th)\1n-\1n\bar X(\th)])\b d\b d\a\){\bf1}_{[\bar\t_j,T]}(\th)\to B_{xx}(\th),\qq\as\\
\ns\ds\Si^{i,\e,\bar\e}_{xx}(\th)=\sum_{j\ges0}\( \int_0^1\3n\int_0^1\1n \si_{xx}(\bar\t_j,\th,\bar X(\th)\1n+\1n\a\b[X^{i,\e,\bar\e}(\th)\1n-\1n\bar X(\th)])\b d\b d\a\){\bf1}_{[\bar\t_j,T]}(\th)\to\Si_{xx}(\th),\qq\as\\
\ns\ds\qq\qq\qq\qq\qq\qq\qq\qq\hb{as long as }\lim_{\e,\bar\e\to0}X^{i,\e,\bar\e}(\th)=\bar X(\th),\qq\as\ea\ee
The above suggests we introduce \rf{X_1} and \rf{X_2}. Then, for $m\ges1$,
\bel{|X_1|}\ba{ll}
\ns\ds\dbE\int_t^s|X_1^{i,\e,\bar\e}(r)|^{2m}dr\\
\ns\ds=\dbE\[\int_t^s\1n\Big|\1n\int_t^r\3n B_x(\th)X_1^{i,\e,\bar\e}(\th)d\th+\2n\int_t^r\3n\Si_x(\th)X_1^{i,\e,\bar\e}(\th)
\1n-\1n\si(\bar\t_i,\th)
{\bf1}_{[\bar\t_i,\bar\t_i+\bar\e)}(\th)dW(\th)\1n-\1n\bar\xi_i{\bf1}_{[\bar\t_i,\bar\t_i+\bar\e)}(r)
\Big|^{2m}\2n dr\]\\
\ns\ds\les K\dbE\Big\{\int_t^s\(\int_t^r|X_1^{i,\e,\bar\e}(\th)|d\th\)^2dr+\int_t^s\[\int_t^r
\(|X_1^{i,\e,\bar\e}
(\th)|^2+{\bf1}_{[\bar\t_i,\bar\t_i+\bar\e)}(\th)\)d\th\]^mdr+\bar\e\Big\}\\
\ns\ds\les K\[\int_t^s\(\dbE\int_t^r|X_1^{i,\e,\bar\e}(\th)|^{2m}d\th\)+\bar\e\].\ea\ee
By Gronwall's inequality, we have the second in \rf{|X-X|,|X_1|}. On the other hand, we have
$$\ba{ll}
\ns\ds X^{i,\e,\bar\e}(s)-\bar X(s)-X_1^{i,\e,\bar\e}(s)\\
\ns\ds=\int_t^s\(B_x(\th)[X^{i,\e,\bar\e}(\th)-\bar X(\th)-X_1^{i,\e,\bar\e}(\th)]+[X^{i,\e,\bar\e}(\th)-\bar X(\th)]^\top B_{xx}^{i,\e,\bar\e}(\th)[X^{i,\e,\bar\e}(\th)-\bar X(\th)]\)d\th\qq\qq\qq\qq\\
\ns\ds\q+\int_t^s\(\Si_x(\th)[X^{i,\e,\bar\e}(\th)-\bar X(\th)-X_1^{i,\e,\bar\e}(\th)]+[X^{i,\e,\bar\e}(\th)-\bar X(\th)]^\top\Si_{xx}^{i,\e,\bar\e}(\th)
[X^{i,\e,\bar\e}(\th)-\bar X(\th)]\)dW(\th)\\
\ns\ds\q-\int_t^sb(\bar\t_i,\th,X^{i,\e,\bar\e}(\th))
{\bf1}_{[\bar\t_i,\bar\t_i+\bar\e)}(\th)d\th+\bar\e\(\int_t^sb_\t^{i,\e,\bar\e}(\th)d\th+\int_t^s\si_\t^{i,\e,\bar\e}(\th)
dW(\th)\)+\e\xi_i{\bf1}_{[\bar\t_i+\bar\e,T]}(s).\ea$$
This leads to
$$\ba{ll}
\ns\ds\dbE\int_t^s|X^{i,\e,\bar\e}(r)-\bar X(r)-X^{i,\e,\bar\e}_1(r)|^{2m}dr\\
\ns\ds=\dbE\int_t^s\Big|\int_t^r\(B_x(\th)[X^{i,\e,\bar\e}(\th)-\bar X(\th)-X_1^{i,\e,\bar\e}(\th)]d\th+\int_t^r\([X^{i,\e,\bar\e}(\th)-\bar X(\th)]^\top B^{i,\e,\bar\e}_{xx}(\th)[X^{i,\e,\bar\e}(\th)-\bar X(\th)]\)d\th\\
\ns\ds\qq+\int_t^r\Si_x(\th)
[X^{i,\e,\bar\e}(\th)-\bar X(\th)-X_1^{i,\e,\bar\e}(\th)]dW(\th)+\int_t^r\([X^{i,\e,\bar\e}(\th)-\bar X(\th)]^\top\Si_{xx}^{i,\e,\bar\e}(\th)[X^{i,\e,\bar\e}(\th)-\bar X(\th)]\)dW(\th)\ea$$
$$\ba{ll}
\ns\ds\qq-
\int_t^rb(\bar\t_i,\th,X^{i,\e,\bar\e}(\th))
{\bf1}_{[\bar\t_i,\bar\t_i+\bar\e)}(\th)d\th+\bar\e\(\int_t^rb_\t^{i,\e,\bar\e}(\th)d\th+\int_t^s\si_\t^{i,\e,\bar\e}(\th)
dW(\th)\)+\e\xi_i{\bf1}_{[\bar\t_i+\bar\e,T]}(r)|^{2m}dr\ea$$
$$\ba{ll}
\ns\ds\les K\Big\{\dbE\int_t^s\[\(\int_t^r|X^{i,\e,\bar\e}(\th)-\bar X(\th)-X_1^{i,\e,\bar\e}(\th)|d\th\)^{2m}+\(\int_t^r|X^{i,\e,\bar\e}(\th)-\bar X(\th)-X_1^{i,\e,\bar\e}(\th)|^2d\th\)^m\\
\ns\ds\qq\qq+\(\int_t^r|X^{i,\e,\bar\e}(\th)-\bar X(\th)|^2d\th\)^{2m}
+\(\int_t^r|X^{i,\e,\bar\e}(\th)-\bar X(\th)|^4d\th\)^m\]dr+\bar\e^{2m}+\e^{2m}\Big\}\\
\ns\ds\les K\[\int_t^s\(\dbE\int_t^r|X^{i,\e,\bar\e}(\th)-\bar X(\th)-X^{i,\e,\bar\e}_1(\th)|^{2m}d\th\)dr+\bar\e^{2m}+\e^{2m}\].\ea$$
Thus, by Gronwall's inequality, we have the first of \rf{|X-X-X_1|,|X_2|}.
Similarly, by the choice of \rf{X_2}, we have
$$\ba{ll}
\ns\ds\dbE\int_t^s|X^{i,\e,\bar\e}_2(r)|^{2m}dr=\int_t^s\dbE\Big|\int_t^r\(B_x
(\th)X^{i,\e,\bar\e}_2(\th)+
X^{i,\e,\bar\e}_1(\th)^\top B_{xx}(\th)X_1^{i,\e,\bar\e}(\th)-b(\bar\t_i,\th){\bf1}_{[\bar\t_i,\bar\t_i+\bar\e)}
(\th)\)d\th\\
\ns\ds\qq\qq+\int_t^r\(\Si_x(\th)X^{i,\e,\bar\e}_2(\th)+X^{i,\e,\bar\e}_1
(\th)^\top\Si_{xx}(\th)
X^{i,\e,\bar\e}_1(\th)\)dW(\th)\\
\ns\ds\qq\qq+\bar\e\(\int_t^sb^i_\t(\th)d\th
+\int_t^s\si^i_\t(\th)dW(\th)\)+\e\xi_i{\bf1}_{[\bar\t_i,T]}(r)\Big|^{2m}dr\\
\ns\ds\qq\les K\dbE\Big\{\int_t^s\[\int_t^r\(|X_2^{i,\e,\bar\e}(\th)|+|X_1^{i,\e,\bar\e}(\th)|^2
+{\bf1}_{[\bar\t_i,\bar\t_i+\bar\e)}(\th)d\th\)^{2m}\\
\ns\ds\qq\qq+\2n\int_t^r\3n\(|X_2^{i,\e,\bar\e}(\th)|^2\1n+\1n|X_1^{i,\e,\bar\e}
(\th)|^4\)^md\th+\bar\e^{2m}+\e^{2m}\Big\}dr\les\1n K\1n\[\int_t^s\3n\(\dbE\2n\int_t^r\3n|X_2^{i,\e,\bar\e}(\th)|^{2m}\1n d\th\)dr\1n+\1n\e^{2m}\1n+\1n\bar\e^{2m}\].\ea$$
By Gronwall's inequality, one has the first of \rf{|X-X-X_1|,|X_2|}. Finally, by \rf{b to} and \rf{B to},

$$\ba{ll}
\ns\ds X^{i,\e,\bar\e}(s)-\bar X(s)-X_1^{i,\e,\bar\e}(s)-X_2^{i,\e,\bar\e}(s)\\
\ns\ds=\int_t^s\(B_x(\th)[X^{i,\e,\bar\e}(\th)-\bar X(\th)-X_1^{i,\e,\bar\e}(\th)-X_2^{i,\e,\bar\e}(\th)]\\
\ns\ds\qq+[X^{i,\e,\bar\e}(\th)-\bar X(\th)]^\top B^{i,\e,\bar\e}_{xx}(\th)[X^{i,\e,\bar\e}(\th)-\bar X(\th)]-X_1^{i,\e,\bar\e}(\th)^\top B_{xx}(\th)X_1^{i,\e,\bar\e}(\th))\)d\th\\
\ns\ds\qq+\int_t^s\(\Si_x(\th)[X^{i,\e,\bar\e}(\th)-\bar X(\th)-X_1^{i,\e,\bar\e}(\th)-X_2^{i,\e,\bar\e}(\th)]\\
\ns\ds\qq+[X^{i,\e,\bar\e}(\th)-\bar X(\th)]^\top\Si^{i,\e,\bar\e}_{xx}(\th)
[X^{i,\e,\bar\e}(\th)-\bar X(\th)]-X_1^{i,\e,\bar\e}(\th)^\top \Si_{xx}(\th)X_1^{i,\e,\bar\e}(\th))\)dW(\th)\\
\ns\ds\qq+\bar\e\[\int_t^s\(b_\t^{i,\e,\bar\e}(\th)-b^i_\t(\th)\)d\th
+\int_t^s\(\si_\t^{i,\e,\bar\e}(\th)
-\si^i_\t(\th)\)dW(\th)\]\\

\ea$$
This leads to the following.
$$\ba{ll}
\ns\ds\dbE\int_t^s|X^{i,\e,\bar\e}(r)-\bar X(r)-X_1^{i,\e,\bar\e}(r)-X_2^{i,\e,\bar\e}(r)|^{2m}dr\\
\ns\ds=\dbE\int_t^s\Big|\int_t^r\(B_x(\th)[X^{i,\e,\bar\e}(\th)-\bar X(\th)-X_1^{i,\e,\bar\e}(\th)-X_2^{i,\e,\bar\e}(\th)]\\
\ns\ds\qq+[X^{i,\e,\bar\e}(\th)-\bar X(\th)-X_1^{i,\e,\bar\e}(\th)]^\top B_{xx}(\th)[X^{i,\e,\bar\e}(\th)-\bar X(\th)+X_1^{i,\e,\bar\e}(\th)]\\
\ns\ds\qq+X_1^{i,\e,\bar\e}(\th)[B^{i,\e,\bar\e}_{xx}(\th)-B_{xx}(\th)]
X_1^{i,\e,\bar\e}(\th))\)d\th+\int_t^s\(\Si_x(\th)[X^{i,\e,\bar\e}(\th)-\bar X(\th)-X_1^{i,\e,\bar\e}(\th)-X_2^{i,\e,\bar\e}(\th)]\\
\ns\ds\qq+[X^{i,\e,\bar\e}(\th)-\bar X(\th)-X_1^{i,\e,\bar\e}(\th)]^\top\Si^{i,\e,\bar\e}_{xx}(\th)
[X^{i,\e,\bar\e}(\th)-\bar X(\th)+X_1^{i,\e,\bar\e}(\th)]\\
\ns\ds\qq+X_1^{i,\e,\bar\e}(\th)[\Si^{i,\e,\bar\e}_{xx}(\th)-\Si_{xx}(\th)]
X_1^{i,\e,\bar\e}(\th))\)dW(\th)+\bar\e\(\int_t^sb^i_\t(\th)d\th+\int_t^s
\si^i_\t(\th)dW(\th)\)\\
\ns\ds\qq+\bar\e\[\int_t^s\(b_\t^{i,\e,\bar\e}(\th)-b_\t^i(\th)\)
d\th+\int_t^s\(\si_\t^{i,\e,\bar\e}(\th)
-\si_\t^i(\th)\)dW(\th)\]\Big|^{2m}dr\ea$$
$$\ba{ll}
\ns\ds\les K\Big\{\int_t^s\dbE\[\(\int_t^r|X^{i,\e,\bar\e}(\th)-\bar X(\th)-X_1^{i,\e,\bar\e}(\th)-X_2^{i,\e,\bar\e}(\th)|d\th\)^{2m}\\
\ns\ds\qq+\(\int_t^r\big[|X^{i,\e,\bar\e}(\th)-\bar X(\th)-X_1^{i,\e,\bar\e}(\th)|\,|X^{i,\e,\bar\e}(\th)-\bar X(\th)+X_1^{i,\e,\bar\e}(\th)|\\
\ns\ds\qq+|X_1^{i,\e,\bar\e}(\th)|^2|B_{xx}^{i,\e,\bar\e}(\th)-B_{xx}(\th)|\big]d\th\)^{2m}+\(\int_t^r|X^{i,\e,\bar\e}(\th)-\bar X(\th)-X_1^{i,\e,\bar\e}(\th)-X_2^{i,\e,\bar\e}(\th)|^2d\th\)^m\\
\ns\ds\qq+\(\int_t^r|X^{i,\e,\bar\e}(\th)-\bar X(\th)-X_1^{i,\e,\bar\e}(\th)|^2|X^{i,\e,\bar\e}(\th)-\bar X(\th)+X_1^{i,\e,\bar\e}(\th)|^2d\th\)^m\\
\ns\ds\qq+\(\int_t^r|X_1^{i,\e,\bar\e}(\th)|^4|\Si^{i,\e,\bar\e}_{xx}(\th)
-\Si_{xx}(\th)|^2d\th\)^m+\bar\e^{2m}\\
\ns\ds\qq+\bar\e^{2m}\[\(\int_t^s|b_\t^{i,\e,\bar\e}(\th)-b_\t^i(\th)|d\th\)^{2m}
+\bar\e^{2m}\(\int_t^s|\si_\t^{i,\e,\bar\e}(\th)
-\si_\t^i(\th)|^2d(\th)\)^m\]dr\Big\}.\ea$$
Hence,
$$\ba{ll}
\ns\ds\dbE\int_t^s|X^{i,\e,\bar\e}(r)-\bar X(r)-X_1^{i,\e,\bar\e}(r)-X_2^{i,\e,\bar\e}(r)|^{2m}dr\\
\ns\ds\les K\dbE\Big\{\int_t^s\(\dbE\int_t^r|X^{i,\e,\bar\e}(\th)-\bar X(\th)-X_1^{i,\e,\bar\e}(\th)-X_2^{i,\e,\bar\e}(\th)|^{2m}d\th\)\\
\ns\ds\q+O(\e^{2m}+\bar\e^{2m})\big[O(\e^{2m}+\bar\e)+O(\bar\e)\big]+O(\bar\e^{2m})
\[\(\int_t^T|B_{xx}^{i,\e,\bar\e}(\th)-B_{xx}(\th)|^{2m\over2m-1}d\th\)\\
\ns\ds\q+\(\int_t^r\3n|\Si^{i,\e,\bar\e}_{xx}(\th)-\Si_{xx}(\th)|^{2m\over m-1}d\th\)\]+\bar\e^{2m}\[\int_t^s\3n|b^{i,\e,\bar\e}_\t(\th)-b_\t^i(\th)|
d\th\)^{2m}\3n
+\(\int_t^s\3n|\si_t^{i,\e,\bar\e}(\th)-\si_t^i(\th)|^2d\th\)^m\]dr\Big\}\\
\ns\ds\les K\[\int_t^s\(\dbE\int_t^r|X^{i,\e,\bar\e}(\th)-\bar X(\th)-X_1^{i,\e,\bar\e}(\th)-X_2^{i,\e,\bar\e}(\th)|^{2m}d\th\)
+o(\e^{2m}+\bar\e^{2m})\].\ea$$
By Gronwall's inequalty, we have \rf{|X-X-X-X|}. \endpf

\ms

\subsection{Variational inequality}

Our next result is the following variational inequality.

\bp{4.4} \sl Let {\rm(H1)--(H2)} and {\rm(H4)} hold. Let $(\t_0,t,x)\in\D_*[0,T]\times\dbR^n$ be given. If $\bar\xi(\cd)\in\BBXi[t,T]$ is an optimal impulse control of form \rf{bar xi} with the perturbation \rf{xi^e}. Then the following holds:
\bel{ineq1}\ba{ll}
\ns\ds0\les\dbE\[H_x[X^{i,\e,\bar\e}_1(T)+X_2^{i,\e,\bar\e}(T)]
+X^{i,\e,\bar\e}_1(T)^\top H_{xx}X^{i,\e,\bar\e}_1(T)\\
\ns\ds\qq+\2n\int_t^s\3n\( G_x(r)[X^{i,\e,\bar\e}_1(r)+X_2^{i,\e,\bar\e}(r)]+X^{i,\e,\bar\e}_1(r)
^\top G_{xx}(r)X^{i,\e,\bar\e}_1(r)\)dr\\
\ns\ds\qq+\bar\e\(\ell_\t(\bar\t_i,\bar\xi_i)+\int_t^Tg^i_\t(r)dr\)
+\e\ell_\xi(\bar\t_i,\bar\xi_i)\xi_i-\int_t^Tg(\bar\t_i,r,\bar X(r)){\bf1}_{[\bar\t_i,\bar\t_i+\bar\e)}(r)dr\]+o(\e+\bar\e),\ea\ee
where
\bel{G}\left\{\2n\ba{ll}
\ns\ds H_x=h_x(\bar X(T)),\qq H_{xx}={1\over2}h_{xx}(\bar X(T)),\\
\ns\ds g^i_\t(r)=g_\t(\bar\t_i,r,\bar X(r)){\bf1}_{[\bar\t_i,T]},\qq i\ges1,\\ [2mm]
\ns\ds G_x(r)=\sum_{j\ges0}g_x(\bar\t_j,r,\bar X(r)){\bf1}_{[\bar\t_j,T]}(r),\q
G_{xx}(r)={1\over2}\sum_{j\ges0}g_{xx}(\bar\t_j,r,\bar X(r)){\bf1}_{[\bar\t_j,T]}(r).\ea\right.\ee

\ep

\it Proof. \rm By the optimality of $\bar\xi(\cd)$, we have
$$\ba{ll}
\ns\ds0\les J(\t_0,t,x;\xi^{i,\e,\bar\e}(\cd))-J(\t_0,t,x;\bar\xi(\cd))\\
\ns\ds\q=\dbE\[h(X^{i,\e,\bar\e}(T))-h(\bar X(T))+\(\ell(\bar\t_i+\bar\e,\bar\xi_i+\e\xi_i)-\ell(\bar\t_i,
\bar\xi_i)\)+\sum_{j\ges0,j\ne i}\int_t^T\(g(\bar\t_j,r,X^{i,\e,\bar\e}(r))
{\bf1}_{[\bar\t_j,T]}(r)\\
\ns\ds\qq\qq\qq\qq\qq\qq+g(\bar\t_i+\bar\e,r,X^{i,\e,\bar\e}(r)){\bf1}_{[\bar\t_i+\bar\e,T]}(r)
\)dr-\sum_{j\ges0}\int_t^Tg(\bar\t_j,r,\bar X(r)){\bf1}_{[\bar\t_j,T]}
(r)dr\]\\
\ns\ds\q=\dbE\Big\{h_x(\bar X(T))[X^{i,\e,\bar\e}(T)-\bar X(T)]\\
\ns\ds\qq+[X^{i,\e,\bar\e}(T)-\bar X(T)]^\top\(\int_0^1\3n\int_0^1h_{xx}(\bar X(T)+\b\a[X^{i,\e,\bar\e}(T)-\bar X(T)])\b d\b d\a\)[X^{i,\e,\bar\e}(T)-\bar X(T)]\\
\ns\ds\qq+\(\int_0^1\ell_\t(\bar\t_i+\a\bar\e,\bar\xi_i+\a\e\xi_i)d\a\)\bar\e
+\(\int_0^1\ell_\xi(\bar\t_i+\a\bar\e,\bar\xi_i+\a\e\xi_i)d\a\)\e\xi_i\\
\ns\ds\qq+\int_t^T\sum_{j\ges0}\(g(\bar\t_j,r,X^{i,\e,\bar\e}(r))
-g(\bar\t_j,r,\bar X(r))\){\bf1}_{[\bar\t_j,T]}(r)dr\\
\ns\ds\qq+\2n\int_t^T\3n\[\(g(\bar\t_i+\bar\e,r,X^{i,\e,\bar\e}(r))
-g(\bar\t_i,r,\bar X(r))\){\bf1}_{[\bar\t_i+\bar\e,
T]}(r)-g(\bar\t_i,r,\bar X(r)){\bf1}_{[\bar\t_i,\bar\t_i+\bar\e)}(r)\]dr\Big\}\\
\ns\ds\q=\dbE\Big\{H_x[X^{i,\e,\bar\e}(T)-\bar X(T)]+[X^{i,\e,\bar\e}(T)-\bar X(T)]^\top H_{xx}^{i,\e,\bar\e}[X^{i,\e,\bar\e}(T)-\bar X(T)]\\
\ns\ds\qq+\ell_\t^{i,\e,\bar\e}\bar\e+\ell_\xi^{i,\e,\bar\e}\e\xi_i
+\int_t^T\sum_{j\ges0}\[g_x(\bar\t_j,r,\bar X(r))[X^{i,\e,\bar\e}(r)-\bar X(r)]\\
\ns\ds\qq+[X^{i,\e,\bar\e}(r)\1n-\1n\bar X(r)]^\top\1n\(\2n\int_0^1\3n\int_0^1 g_{xx}(\bar\t_j,r,\bar X(r)\1n+\1n\a\b[X^{i,\e,\bar\e}(r)\1n-\1n\bar X(r)])\b d\a d\b\)[X^{i,\e,\bar\e}(r)\1n-\1n\bar X(r)]\]{\bf1}_{[\bar\t_j,T]}(r)dr\\
\ns\ds\qq+\bar\e\(\2n\int_t^T\3n\int_0^1\3n g_\t(\bar\t_i+\a\bar\e,r,\bar X(r))d\a\){\bf1}_{[\bar\t_i+\bar\e,T]}(r)dr-\int_t^Tg(\bar\t_i,r,\bar X(r)){\bf1}_{[\bar\t_i,\bar\t_i+\bar\e)}(r)dr\Big\}\\
\ns\ds\q=\dbE\Big\{H_x[X^{i,\e,\bar\e}(T)-\bar X(T)]+[X^{i,\e,\bar\e}(T)-\bar X(T)]^\top H_{xx}^{i,\e,\bar\e}[X^{i,\e,\bar\e}(T)-\bar X(T)]\\
\ns\ds\qq+\ell_\t^{i,\e,\bar\e}\bar\e+\ell_\xi^{i,\e,\bar\e}\e\xi_i
+\int_t^TG_x(r)[X^{i,\e,\bar\e}(r)-\bar X(r)]+[X^{i,\e,\bar\e}(r)-\bar X(r)]^\top G_{xx}^{i,\e,\bar\e}(r)[X^{i,\e,\bar\e}(r)-\bar X(r)]\\
\ns\ds\qq+\bar\e\int_t^Tg_\t^{i,\bar\e}(r){\bf1}_{[\bar\t_i+\bar\e,T]}(r)dr
-\int_t^Tg(\bar\t_i,r,\bar X(r)){\bf1}_{[\bar\t_i,\bar\t_i+\bar\e)}(r)dr\Big\}+o(\e+\bar\e),\ea$$
where
$$\left\{\2n\ba{ll}
\ns\ds H^{i,\e,\bar\e}_{xx}=\int_0^1\int_0^1h_{xx}(\bar X(T)+\b\a(X^{i,\e,\bar\e}(T)-\bar X(T])\b d\b d\a\to H_{xx},\\
\ns\ds\ell^{i,\e,\bar\e}_\t=\int_0^1\ell_\t(\bar\t_i+\a\bar\e,\bar\xi_i+\a\e\xi_i)
d\a\to\ell_\t(\bar\t_i,\bar\xi_i),
\q\ell^{i,\e,\bar\e}_\xi=\int_0^1\ell_\xi(\bar\t_i+\a\bar\e,\bar\xi_i+\a\e\xi_i)
d\a\to\ell_\xi(\bar\t_i,\bar\xi_i),\\
\ns\ds g_\t^{i,\bar\e}(r)\to g_\t^i(r),\qq G^{i,\e,\bar\e}_{xx}(r)=\int_0^1\int_0^1g_{xx}(\bar\t_j,r,\bar X(r)+\b\a[X^{i,\e,\bar\e}(r)-\bar X(r)]\b d\b d\a\to G_{xx}(r).\ea\right.$$
Thus, we have \rf{ineq1} by the estimate \rf{|X-X|,|X_1|} and \rf{|X-X-X_1|,|X_2|}. \endpf

\ms

We note that the variational inequality \rf{ineq1} is not very useful since
it involves $X^{i,\e,\bar\e}_1(\cd)$ and $X^{i,\e,\bar\e}_2(\cd)$, which are
depending on the variation $\xi^{i,\e,\bar\e}(\cd)$ of the optimal impulse control and they might not even be continuous. Thus, we need to go one more step.

\subsection{Duality principles}

In order to introduce the duality principles to eliminate the dependence on the $X_1^{i,\e,\bar\e}(\cd)$ and $X_2^{i,\e,\bar\e}(\cd)$ in \rf{ineq1}, we now define (compare with \rf{X_1} and \rf{X_2})
\bel{h X_1}\ba{ll}
\ns\ds\h X_1^{i,\e,\bar\e}(s)=\int_t^sB_x(r) X^{i,\e,\bar\e}_1(r)dr+\int_t^s\(\Si_x(r)
X^{i,\e,\bar\e}_1(r)-\si(\bar\t_i,r){\bf1}_{[\bar\t_i,
\bar\t_i+\bar\e)}(r)\)dW(r)\\
\ns\ds\qq\qq\equiv X_1^{i,\e,\bar\e}(s)+\bar\xi_i{\bf1}_{[\bar\t_i,\bar\t_i+\bar\e)}(s),\qq s\in[t,T],\ea\ee
and
\bel{h X_2}\ba{ll}
\ns\ds\h X^{i,\e,\bar\e}_2(s)=\int_t^s\(B_x(r)X^{i,\e,\bar\e}_2(r)+X^{i,\e,\bar\e}_1(r)^\top B_{xx}(r)X_1^{i,\e,\bar\e}(r)-b(\bar\t_i,r){\bf1}_{[\bar\t_i,\bar\t_i+\bar\e)}
(r)\)dr\\
\ns\ds\qq\qq\qq+\int_t^s\(\Si_x(r)X^{i,\e,\bar\e}_2(r)+X^{i,\e,\bar\e}_1
(r)^\top\Si_{xx}(r)X^{i,\e,\bar\e}_1(r)\)dW(r)\\
\ns\ds\qq\qq\equiv X_2^{i,\e,\bar\e}(s)-\bar\e\z^i(s)-\e\xi_i{\bf1}_{[\bar\t_i+\bar\e,T]}(s),\q s\in[t,T],\ea\ee
They are It\^o processes, by which we mean that they can be written as the sums of Lebesgue integral (the drift term) and the It\^o's integral (diffusion term). Therefore, they are continuous in $s$. Clearly,
\bel{X_1+}X_1^{i,\e,\bar\e}(s)=\h X_1^{i,\e,\bar\e}(s)-\bar\xi_i{\bf1}_{[\bar\t_i,\bar\t_i+\bar\e)}(s),\ee
and
\bel{X_2+}X_2^{i,\e,\bar\e}(s)=\h X_2^{i,\e,\bar\e}(s)+\bar\e\z^i(s)
+\e\xi_i{\bf1}_{[\bar\t_i+\bar\e,T]}(s).\ee
Note that since on $[t,\bar\t_i+\bar\e)$, the optimal impulse control
is not perturbed yet (the perturbation is going to be made at $\bar\t_i+\bar\e$ from $\bar\xi_i$ to $\bar\xi_i+\e\xi_i$), which leads to
$$0=\bar\xi_i^{\,\top}{\bf1}_{[\bar\t_i,\bar\t_i+\bar\e)}(r)G_{xx}(r) X_1^{i,\e,\bar\e}(r)=\bar\xi_i^{\,\top}{\bf1}_{[\bar\t_i,\bar\t_i+\bar\e)}(r)
G_{xx}(r)[\h X_1^{i,\e,\bar\e}(r)-\bar\xi_i{\bf1}_{[\bar\t_i,\bar\t_i+\bar\e)}(s)].$$
Thus,
$$\bar\xi_i^{\,\top}{\bf1}_{[\bar\t_i,\bar\t_i+\bar\e)}(r)
G_{xx}(r)\h X_1^{i,\e,\bar\e}(r)=\bar\xi_i^{\,\top}G_{xx}(r)\bar\xi_i
{\bf1}_{[\bar\t_i,\bar\t_i+\bar\e)}(s).$$
Also, since we have assumed $\bar\t_{\k(\bar\xi(\cd))}<T$ (the other case can be treated similarly), one must have
$$\bar\t_i+\bar\e<T,\qq 1\les i\les\k(\bar\xi(\cd)),$$
if $\bar\e$ is small. Consequently,
$${\bf1}_{[\bar\t_i,\bar\t_i+\bar\e)}(T)=0,\qq1\les i\les\k(\bar\xi(\cd)).$$
Hence, the variational inequality \rf{ineq1} reads
\bel{ineq2}\ba{ll}
\ns\ds0\les\dbE\[H_x[\h X^{i,\e,\bar\e}_1(T)-\bar\xi_i{\bf1}_{[\bar\t_i,\bar\t_i+\bar\e)}(T)
+\h X_2^{i,\e,\bar\e}(T)+\bar\e\z^i(T)+\e\xi_i{\bf1}_{[\bar\t_i,T]}(T)]\\
\ns\ds\qq+[\h X_1^{i,\e,\bar\e}(T)-\bar\xi_i{\bf1}_{[\bar\t_i,\bar\t_i+\bar\e)}(T)]^\top H_{xx}[\h X_1^{i,\e,\bar\e}(T)-\bar\xi_i{\bf1}_{[\bar\t_i,\bar\t_i+\bar\e)}(T)]\\
\ns\ds\qq+\2n\int_t^T\3n G_x(r)[\h X^{i,\e,\bar\e}_1(r)-\bar\xi_i{\bf1}_{[\bar\t_i,\bar\t_i+\bar\e)}(r)+\h X_2^{i,\e,\bar\e}(r)+\bar\e\z^i(r)+\e\xi_i{\bf1}_{[\bar\t_i,T]}(r)]\\
\ns\ds\qq+[\h X^{i,\e,\bar\e}_1(r)-\bar\xi_i{\bf1}_{[\bar\t_i,\bar\t_i+\bar\e)}(r)
]^\top G_{xx}(r)[\h X^{i,\e,\bar\e}_1(r)-\bar\xi_i{\bf1}_{[\bar\t_i,\bar\t_i+\bar\e)}(r)]dr\]\\
\ns\ds\qq+\bar\e\(\ell_\t(\bar\t_i,\bar\xi_i)+\int_t^Tg^i_\t(r)dr\)
+\e\ell_\xi(\bar\t_i,\bar\xi_i)\xi_i-\int_t^Tg(\bar\t_i,r,\bar X(r)){\bf1}_{[\bar\t_i,\bar\t_i+\bar\e)}(r)dr+o(\e+\bar\e)\\
\ns\ds=\dbE\[H_x[\h X^{i,\e,\bar\e}_1(T)+\h X_2^{i,\e,\bar\e}(T)]+\h X^{i,\e,\bar\e}_1(T)^\top H_{xx}\h X^{i,\e,\bar\e}_1(T)+\bar\e H_x\z^i(T)+\e  H_x\xi_i\\
\ns\ds\qq+\2n\int_t^T\(G_x(r)[\h X^{i,\e,\bar\e}_1(r)+\h X_2^{i,\e,\bar\e}(r)]+\h X^{i,\e,\bar\e}_1(r)^\top G_{xx}(r)\h X^{i,\e,\bar\e}_1(r)-g(\bar\t_i,r,\bar X(r)){\bf1}_{[\bar\t_i,\bar\t_i+\bar\e)}(r)\\
\ns\ds\qq-G_x(r)\bar\xi_i{\bf1}_{[\bar\t_i,\bar\t_i+\bar\e)}(r)+\bar\e G_x(r)\z^i(r)+\e G_x(r)\xi_i{\bf1}_{[\bar\t_i,T]}(r)-\bar
\xi_i^{\,\top}G_{xx}(r)\bar\xi_i{\bf1}_{[\bar\t_i,\bar\t_i+\bar\e)}(r)\\
\ns\ds\qq+\bar\e\(\ell_\t(\bar\t_i,\bar\xi_i)+\int_t^Tg^i_\t(r)dr\)
+\e\ell_\xi(\bar\t_i,\bar\xi_i)\xi_i\]+o(\e+\bar\e),\ea\ee
which leads to 
\bel{ineq3}\ba{ll}
\ns\ds0\les\dbE\Big\{H_x[\h X^{i,\e,\bar\e}_1(T)+\h X_2^{i,\e,\bar\e}(T)]+\h X^{i,\e,\bar\e}_1(T)^\top H_{xx}\h X^{i,\e,\bar\e}_1(T)\\
\ns\ds\qq+\2n\int_t^T\(G_x(r)[\h X^{i,\e,\bar\e}_1(r)+\h X_2^{i,\e,\bar\e}(r)]+\h X^{i,\e,\bar\e}_1(r)^\top G_{xx}(r)\h X^{i,\e,\bar\e}_1(r)\)dr\\
\ns\ds\qq-\int_t^T\(G_x(r)\bar\xi_i+\bar
\xi_i^{\,\top}G_{xx}(r)\bar\xi_i+g(\bar\t_i,r,\bar X(r))\){\bf1}_{[\bar\t_i,\bar\t_i+\bar\e)}(r)dr\\
\ns\ds\qq+\bar\e\(\ell_\t(\bar\t_i,\bar\xi_i)+H_x\z^i(T)+\int_t^T[G_x(r)\z^i(r)+
g^i_\t(r)]dr\)\\
\ns\ds\qq+\e\(H_x+\ell_\xi(\bar\t_i,\bar\xi_i)+\int_t^TG_x(r)
{\bf1}_{[\bar\t_i,T]}(r)dr\)\xi_i\Big\}+o(\e+\bar\e),\ea\ee
We note that the above still depends on $\xi^{i,\e,\bar\xi}(\cd)$, but through two It\^o processes. We now want to use duality principles to get rid of such dependence. To this end, by the definitions, we obtain the SDEs for
$\h X_1^{i,\e,\bar\e}(\cd)$ and $\h X_2^{i,\e,\bar\e}(\cd)$:
\bel{h X_1*}\ba{ll}
\ns\ds\h X_1^{i,\e,\bar\e}(s)=\int_t^s\(B_x(r)\h X^{i,\e,\bar\e}_1(r)-B_x(r)\bar\xi_i{\bf1}_{[\bar\t_i,\bar\t_i+\bar\e)}(r)\)dr\\
\ns\ds\qq\qq\qq+\int_t^s\(\Si_x(r)\h
X^{i,\e,\bar\e}_1(r)-[\Si_x(r)\bar\xi_i+\si(\bar\t_i,r)]{\bf1}_{[\bar\t_i,
\bar\t_i+\bar\e)}(r)\)dW(r)\\
\ns\ds\qq\qq\equiv\int_t^s\(B_x(r)\h X^{i,\e,\bar\e}_1(r)+\h b_{1,\bar\e}^i(r)\)dr+\int_t^s\(\Si_x(r)\h
X^{i,\e,\bar\e}_1(r)+\h\si_{1,\bar\e}^i(r)\)dW(r),\ea\ee
with
\bel{b_1}\h b_{1,\bar\e}^i(r)=-B_x(r)\bar\xi_i{\bf1}_{[\bar\t_i,\bar\t_i+\bar\e)}(r),\qq
\h\si_{1,\bar\e}^i(r)=-[\Si_x(r)\bar\xi_i+\si(\bar\t_i,r)]{\bf1}_{[\bar\t_i,
\bar\t_i+\bar\e)}(r),\ee
and
$$\ba{ll}
\ns\ds\h X^{i,\e,\bar\e}_2(s)=\int_t^s\(B_x(r)X^{i,\e,\bar\e}_2(r)+X^{i,\e,\bar\e}_1(r)^\top B_{xx}(r)X_1^{i,\e,\bar\e}(r)-b(\bar\t_i,r){\bf1}_{[\bar\t_i,\bar\t_i+\bar\e)}
(r)\)dr\\
\ns\ds\qq\qq\qq+\int_t^s\(\Si_x(r)X^{i,\e,\bar\e}_2(r)+X^{i,\e,\bar\e}_1
(r)^\top\Si_{xx}(r)
X^{i,\e,\bar\e}_1(r)\)dW(r)\\
\ns\ds\qq\qq=\int_t^s\(B_x(r)[\h X^{i,\e,\bar\e}_2(r)+\bar\e\z^i(r)+\e\xi_i{\bf1}_{[\bar\t_i+\bar\e)}(r)]
-b(\bar\t_i,r){\bf1}_{[\bar\t_i,\bar\t_i+\bar\e)}
(r)\\
\ns\ds\qq\qq\qq+[\h X^{i,\e,\bar\e}_1(r)-\bar\xi_i{\bf1}_{[\bar\t_i,\bar\t_i+\bar\e)}(r)]^\top B_{xx}(r)[\h X_1^{i,\e,\bar\e}(r)-\bar\xi_i{\bf1}_{[\bar\t_i,\bar\t_i+\bar\e)}(r)]\)dr\\
\ns\ds\qq\qq\qq+\int_t^s\(\Si_x(r)[\h X^{i,\e,\bar\e}_2(r)+\bar\e\z^i(r)+\e\xi_i{\bf1}_{[\bar\t_i,\bar\t_i+\bar\e)}(r)]\\
\ns\ds\qq\qq\qq+[\h X_1^{i,\e,\bar\e}(r)-\bar\xi_i{\bf1}_{[\bar\t_i,\bar\t_i+\bar\e)}(r)]^\top
\Si_{xx}(r)[\h X_1^{i,\e,\bar\e}(r)-\bar\xi_i{\bf1}_{[\bar\t_i,\bar\t_i+\bar\e)}(r)]\)dW(r)\\
\ns\ds\qq\qq=\int_t^s\(B_x(r)\h X^{i,\e,\bar\e}_2(r)+\bar\e B_x(r)\z^i(r)+\e B_x(r)\xi_i{\bf1}_{[\bar\t_i,\bar\t_i+\bar\e)}(r)
-b(\bar\t_i,r){\bf1}_{[\bar\t_i,\bar\t_i+\bar\e)}
(r)\\
\ns\ds\qq\qq\qq+\h X^{i,\e,\bar\e}_1(r)^\top B_{xx}(r)\h X_1^{i,\e,\bar\e}(r)-\bar\xi_i^{\,\top}B_{xx}(r)\bar\xi_i{\bf1}_{[\bar\t_i,\bar\t_i+\bar\e)}(r)]\)dr\\
\ns\ds\qq\qq\qq+\int_t^s\(\Si_x(r)\h X^{i,\e,\bar\e}_2(r)+\bar\e\Si_x(r)\z^i(r)+\e\Si_x(r)\xi_i{\bf1}_{[\bar\t_i,\bar\t_i+\bar\e)}(r)
\\
\ns\ds\qq\qq\qq+\h X_1^{i,\e,\bar\e}(r)^\top
\Si_{xx}(r)\h X_1^{i,\e,\bar\e}(r)-\bar\xi_i^\top\Si_{xx}(r)\bar\xi_i
{\bf1}_{[\bar\t_i,\bar\t_i+\bar\e)}(r)]\)dW(r),\ea$$
which gives
\bel{X_2*}\h X_2^{i,\e,\bar\e}(s)=\int_t^s\(B_x(r)\h X^{i,\e,\bar\e}_2(r)+\h b_{2,\e,\bar\e}^i(r)\)dr+\int_t^s\(\Si_x(r)\h X^{i,\e,\bar\e}_2(r)+\h\si_{2,
\e,\bar\e}^i(r)\)dW(r),\ee
with
\bel{b_2}\ba{ll}
\ns\ds\h b_{2,\e,\bar\e}^i(r)=\bar\e B_x(r)\z^i(r)\1n+\1n
\[\e B_x(r)\xi_i\1n-\1n\bar\xi_i^{\,\top}B_{xx}(r)\bar\xi_i\1n-\1n b(\bar\t_i,r)\]{\bf1}_{[\bar\t_i,\bar\t_i
\1n+\1n\bar\e)}(r)\1n+\1n\h X^{i,\e,\bar\e}_1(r)^\top B_{xx}(r)\h X_1^{i,\e,\bar\e}(r),\\
\ns\ds\h\si_{2,\e,\bar\e}^i(r)=\bar\e\Si_x(r)\z^i(r)
+\[\e\Si_x(r)\xi_i-\bar\xi_i^{\,\top}\Si_{xx}(r)\bar\xi_i\]{\bf1}_{[\bar\t_i,
\bar\t_i+\bar\e)}(r)+\h X_1^{i,\e,\bar\e}(r)^\top\Si_{xx}(r)\h X_1^{i,\e,\bar\e}(r).\ea\ee
The following duality relation can be used to get rid of the linear terms in $\h X_1^{i,\e,\bar\e}(\cd)+\h X_2^{i,\e,\bar\e}(\cd)$ on the right-hand side of \rf{ineq3}.

\bp{} \sl Let {\rm(H1)--(H2)} and {\rm(H4)} hold and $(Y(\cd),Z(\cd))$ is the adapted solution of the following BSDE:
\bel{BSDE1}\left\{\2n\ba{ll}
\ds dY(s)=-\[B_x(s)^\top Y(s)+\Si_x(s)^\top Z(s)+G_x(s)^\top\]ds+Z(s)dW(s),\q s\in[t,T],\\
\ns\ds Y(T)=H_x^\top,\ea\right.\ee
where we recall \rf{G} (for the data). Then
\bel{duality1}\ba{ll}
\ns\ds\dbE\[H\(\h X_1^{i,\e,\bar\e}(T)+\h X_2^{i,\e,\bar\e}(T)\)+
\int_t^TG_x(r)\(\h X^{i,\e,\bar\e}_1(r)+\h X_2^{i,\e,\bar\e}(r)\)dr\]\\
\ns\ds=\dbE\int_t^T\(\lan Y(r),\h b_{1,\bar\e}^i(r)+\h b_{2,\e,\bar\e}^i(r)\ran+\lan Z(r),\h\si_{1,\bar\e}^i(r)+\h\si_{2,\e,\bar\e}^i(r)\ran\)dr.\ea\ee

\ep

The key point is that the BSDE \rf{BSDE1} is only depending on the optimal impulse control, and is independent of the variation $\xi^{i,\e,\bar\e}(\cd)$.

\ms

\it Proof. \rm By It\^o's formula, together with some approximate arguments, we have
$$\ba{ll}
\ns\ds\dbE\[\lan Y(s),\h X^{i,\e,\bar\e}_1(s)\ran\]=\dbE\int_t^s\(\lan-B_x(r)^\top Y(r)-\Si_x(r)^\top Z(r)-G_x(r)^\top,\h X_1^{i,\e,\bar\e}(r)\ran\\
\ns\ds\qq\qq\qq\qq\qq\qq+\lan Y(r),B_x(r)\h X_1^{i,\e,\bar\e}(r)+\h b_{1,\bar\e}^i(r)\ran+\lan Z(r),\Si_x(r)X^{i,\e,\bar\e}_1(r)+\h\si_{1,\bar\e}^i(r)\ran\)dr\\
\ns\ds\qq\qq\qq\qq~=\dbE\int_t^s\(-G_x(r)X_1^{i,\e,\bar\e}(r)+\lan Y(r),\h b_{1,\bar\e}^i(r)\ran+\lan Z(r),\h\si_{1,\bar\e}^i(r)\ran\)dr.\ea$$
Then
$$\dbE\[H_x\h X_1^{i,\e,\bar\e}(T)+\int_t^TG_x(r)\h X_1^{i,\e,\bar\e}(r)dr\]=\dbE\int_t^T\(\lan Y(r),\h b_{1,\bar\e}^i(r)\ran+\lan Z(r),\h\si_{1,\bar\e}^i(r)\ran\)dr.$$
Similarly,
$$\ba{ll}
\ns\ds\dbE\[\lan Y(s),\h X_2^{i,\e,\bar\e}(s)\ran\]=\dbE\int_t^s\(\lan-B_x(r)^\top Y(r)-\Si_x(r)^\top Z(r)-G_x(r)^\top,\h X^{i,\e,\bar\e}_2(r)\ran\\
\ns\ds\qq\qq\qq\qq\qq+\lan Y(r),B_x(r)\h X_2^{i,\e,\bar\e}(r)+\h b_{2,\e,\bar\e}^i(r)\ran+\lan Z(r),\Si_x(r)\h X^{i,\e,\bar\e}_2(r)+\h\si^i_{2,\e,\bar\e}(r)\ran\)dr\\
\ns\ds\qq\qq\qq\qq~=\dbE\int_t^s\(-G_x(r)X_2^{i,\e,\bar\e}(r)+\lan Y(r),\h b_{2,
\e,\bar\e}^i(r)\ran+\lan Z(r),\h\si_{2,\e,\bar\e}^i(r)\ran\)dr.\ea$$
Then
$$\dbE\[H_xX_2^{i,\e,\bar\e}(T)+\int_t^TG_x(r)X_2^{i,\e,\bar\e}(r)dr\]=\dbE\int_t^T\(\lan Y(r),\h b_{2,\e,\bar\e}^i(r)\ran+\lan Z(r),\h\si_{2,\e,\e\bar\e}^i(r)\ran\)dr.$$
Consequently, we have \rf{duality1}. \endpf

\ms

Note that
$$\ba{ll}
\ns\ds\dbE\int_t^T\(\lan Y(r),\h b^i_{1,\bar\e}(r)+\h b^i_{2,\e,\bar\e}(r)\ran+\lan Z(r),\h\si^i_{1,\bar\e}(r)+\h\si^i_{2,\e,\bar\e}(r)\ran\)dr\\
\ns\ds=\dbE\int_t^T\(\lan Y(r),-B_x(r)\bar\xi_i{\bf1}_{[\bar\t_i,\bar\t_i+\bar\e)}(r)+\bar\e B_x(r)\z^i(r)+
\big[\e B_x(r)\xi_i-\bar\xi_i^{\,\top}\h B_x(r)\bar\xi_i-b(\bar\t_i,r)\big]{\bf1}_{[\bar\t_i,\bar\t_i+\bar\e)}(r)\ran\\
\ns\ds\qq+\lan\h X^{i,\e,\bar\e}_1(r)Y(r)B_{xx}(r)\h X_1^{i,\e,\bar\e}(r)\ran+\lan Z(r),-\big[\Si_x(r)\bar\xi_i+\si(\bar\t_i,r)\big]{\bf1}_{[\bar\t_i,
\bar\t_i+\bar\e)}(r)\ran\\
\ns\ds\qq+\lan Z(r),\bar\e\Si_x(r)\z^i(r)
+\big[\e\Si_x(r)\xi_i-\bar\xi_i^{\,\top}\Si_{xx}(r)\bar\xi_i\big]{\bf1}_{[\bar\t_i,
\bar\t_i+\bar\e)}(r)\ran+\lan\h X_1^{i,\e,\bar\e}(r),Z(r)\Si_{xx}(r)\h X_1^{i,\e,\bar\e}(r)\ran\)dr\\
\ns\ds=\dbE\int_t^T\(-\lan B_x(r)^\top Y(r)+\Si_x(r)^\top Z(r),\bar\xi_i\ran{\bf1}_{[\bar\t_i,\bar\t_i+\bar\e)}(r)+\e\lan B_x(r)^\top Y(r)+\Si_x(r)^\top Z(r),\xi_i\ran{\bf1}_{[\bar\t_i,
\bar\t_i+\bar\e)}(r)\\
\ns\ds\qq-\big(\lan b(\bar\t_i,r),Y(r)\ran+\lan\si(\bar\t_i,r),Z(r)\big)
{\bf1}_{[\bar\t_i,\bar\t_i+\bar\e)}(r)\\
\ns\ds\qq+\bar\e\lan B_x(r)^\top Y(r)+\Si_x(r)^\top Z(r),\z^i(r)\ran+\lan\h X^{i,\e,\bar\e}_1(r),[Y(r)B_{xx}(r)+Z(r)\Si_{xx}(r)]\h X_1^{i,\e,\bar\e}(r)\ran\\
\ns\ds\qq-\lan\bar\xi_i,[Y(r)B_{xx}(r)+Z(r)\Si_{xx}(r)]\bar\xi_i\ran{\bf1}_{[\bar\t_i,
\bar\t_i+\bar\e)}(r)\)dr\\
\ns\ds=\dbE\int_t^T\[\(-\lan B_x(r)^\top Y(r)+\Si_x(r)^\top Z(r),\bar\xi_i\ran-\big(\lan b(\bar\t_i,r),Y(r)\ran+\lan\si(\bar\t_i,r),Z(r)\big)\\
\ns\ds\qq-\lan\bar\xi_i,[Y(r)B_{xx}(r)+Z(r)\Si_{xx}(r)]\bar\xi_i\ran\){\bf1}_{[\bar\t_i,
\bar\t_i+\bar\e)}(r)+\bar\e\lan B_x(r)^\top Y(r)+\Si_x(r)^\top Z(r),\z^i(r)\ran\\
\ns\ds\qq+\lan\h X^{i,\e,\bar\e}_1(r),[Y(r)B_{xx}(r)+Z(r)\Si_{xx}(r)]\h X_1^{i,\e,\bar\e}(r)\ran\]dr+o(\e+\bar\e).\ea$$
Now, we define
$$\BH(\{\t_j\}_{j\ges0},r,x,y,z)=\sum_{j\ges0}\(y^\top b(\t_j,r,x)+z^\top\si(\t_j,r,x)+g(\t_j,r,x)\){\bf1}_{[\t_j,T]}(r),$$
and
$$\ba{ll}
\ns\ds\BH(r)=\BH(\{\bar\t_j\}_{j\ges0},r,\bar X(r),Y(r),Z(r))\\
\ns\ds\qq\equiv\sum_{j\ges0}\(Y(r)^\top b(\t_j,r,\bar X(r))+Z(r)^\top\si(\t_j,r,\bar X(r))+g(\t_j,r,\bar X(r))\){\bf1}_{[\t_j,T]}(r),\\
\ns\ds\BH_x(r)\equiv\sum_{j\ges0}\(Y(r)^\top b_x(\t_j,r,\bar X(r))+Z(r)^\top\si_x(\t_j,r,\bar X(r))+g_x(\t_j,r,\bar X(r))\){\bf1}_{[\t_j,T]}(r),\\
\ns\ds\BH_{xx}(r)\equiv\sum_{j\ges0}\(Y(r)^\top b_{xx}(\t_j,r,\bar X(r))+Z(r)^\top\si_{xx}(\t_j,r,\bar X(r))+g_{xx}(\t_j,r,\bar X(r))\){\bf1}_{[\t_j,T]}(r).\ea$$
Hence, we may write the above more compactly as
\bel{duality1*}\ba{ll}
\ns\ds\dbE\[H\(\h X_1^{i,\e,\bar\e}(T)+\h X_2^{i,\e,\bar\e}(T)\)+
\int_t^TG_x(r)\(\h X^{i,\e,\bar\e}_1(r)+\h X_2^{i,\e,\bar\e}(r)\)dr\]\\
\ns\ds=\dbE\int_t^T\(\lan Y(r),\h b_{1,\bar\e}^i(r)+\h b_{2,\e,\bar\e}^i(r)\ran+\lan Z(r),\h\si_{1,\bar\e}^i(r)+\h\si_{2,\e,\bar\e}^i(r)\ran\)dr\\
\ns\ds=\dbE\int_t^T\[\(-\lan\BH_x(r)-G_x(r),\bar\xi_i\ran-\big(\BH(r)-G(r)\big)
-\lan\bar\xi_i,[\BH_{xx}(r)-G_{xx}(r)]\bar\xi_i\ran\){\bf1}_{[\bar\t_i,
\bar\t_i+\bar\e)}(r)\\
\ns\ds\qq+\bar\e\lan\BH_x(r)-G_x(r),\z^i(r)\ran+\lan\h X^{i,\e,\bar\e}_1(r),\BH_{xx}(r)-G_{xx}(r)]\h X_1^{i,\e,\bar\e}(r)\ran\]dr+o(\e+\bar\e).\ea\ee
Using the above result, we can rewrite the variational inequality \rf{ineq3} as follows.
\bel{ineq3}\ba{ll}
\ns\ds0\les\dbE\Big\{H_x[\h X^{i,\e,\bar\e}_1(T)+\h X_2^{i,\e,\bar\e}(T)]+\h X^{i,\e,\bar\e}_1(T)^\top H_{xx}\h X^{i,\e,\bar\e}_1(T)\\
\ns\ds\qq+\2n\int_t^T\(G_x(r)[\h X^{i,\e,\bar\e}_1(r)+\h X_2^{i,\e,\bar\e}(r)]+\h X^{i,\e,\bar\e}_1(r)^\top G_{xx}(r)\h X^{i,\e,\bar\e}_1(r)\)dr\\
\ns\ds\qq-\int_t^T\(G_x(r)\bar\xi_i+\bar
\xi_i^{\,\top}G_{xx}(r)\bar\xi_i+g(\bar\t_i,r,\bar X(r))\){\bf1}_{[\bar\t_i,\bar\t_i+\bar\e)}(r)dr\\
\ns\ds\qq+\bar\e\(\ell_\t(\bar\t_i,\bar\xi_i)+H_x\z^i(T)+\int_t^T[G_x(r)\z^i(r)+
g^i_\t(r)]dr\)\\
\ns\ds\qq+\e\(H_x+\ell_\xi(\bar\t_i,\bar\xi_i)+\int_t^TG_x(r)
{\bf1}_{[\bar\t_i,T]}(r)dr\)\xi_i\Big\}+o(\e+\bar\e)\\
\ns\ds=\dbE\Big\{\h X^{i,\e,\bar\e}_1(T)^\top H_{xx}\h X^{i,\e,\bar\e}_1(T)+\int_t^T\h X^{i,\e,\bar\e}_1(r)^\top G_{xx}(r)\h X^{i,\e,\bar\e}_1(r)dr\\
\ns\ds\qq+\int_t^T\(\lan Y(r),\h b_{1,\bar\e}^i(r)+\h b_{2,\e,\bar\e}^i(r)\ran+\lan Z(r),\h\si_{1,\bar\e}^i(r)+\h\si_{2,\e,\bar\e}^i(r)\ran\)dr\\
\ns\ds\qq-\int_t^T\(G_x(r)\bar\xi_i+\bar
\xi_i^{\,\top}G_{xx}(r)\bar\xi_i+g(\bar\t_i,r,\bar X(r))\){\bf1}_{[\bar\t_i,\bar\t_i+\bar\e)}(r)dr\\
\ns\ds\qq+\bar\e\(\ell_\t(\bar\t_i,\bar\xi_i)+H_x\z^i(T)+\int_t^T[G_x(r)\z^i(r)+
g^i_\t(r)]dr\)\\
\ns\ds\qq+\e\(H_x+\ell_\xi(\bar\t_i,\bar\xi_i)+\int_t^TG_x(r)
{\bf1}_{[\bar\t_i,T]}(r)dr\)\xi_i\Big\}+o(\e+\bar\e)\\
\ns\ds\q=\dbE\Big\{\h X^{i,\e,\bar\e}_1(T)^\top H_{xx}\h X^{i,\e,\bar\e}_1(T)+\2n\int_t^T\h X^{i,\e,\bar\e}_1(r)^\top\BH_{xx}(r)\h X^{i,\e,\bar\e}_1(r)dr\\
\ns\ds\q+\int_t^T\[-\lan\BH_x(r),\bar\xi_i\ran-\BH(r)\]
{\bf1}_{[\bar\t_i,\bar\t_i+\bar\e)}(r)dr\\
\ns\ds\qq+\bar\e\(\ell_\t(\bar\t_i,\bar\xi_i)+H_x\z^i(T)+\int_t^T[\lan\BH_x(r),\z^i(r)\ran+
g^i_\t(r)]dr\)\\
\ns\ds\qq+\e\(H_x+\ell_\xi(\bar\t_i,\bar\xi_i)+\int_t^TG_x(r)
{\bf1}_{[\bar\t_i,T]}(r)dr\)\xi_i\Big\}+o(\e+\bar\e)\ea\ee
We note that the above is still depending on the variation $\xi^{i,\e,\bar\e}(\cd)$ through $\h X_1^{i,\e,\bar\e}(\cd)$ (quadratically). Now, we want to use the second duality relation to eliminate such a dependance. To this end, we prove the following.

\bp{} \sl Let {\rm(H1)--(H2)} and {\rm(H4)} hold. Let $(\t_0,t,x)\in\D_*[0,T]\times\dbR^n$ be given and $\bar\xi(\cd)$ be an optimal impulse control. Let $(P(\cd),Q(\cd))$ be the adapted solution of the following BSDE:
\bel{BSDE2}\left\{\2n\ba{ll}
\ns\ds dP(s)=-\[P(s)B_x(s)+B_x(s)^\top P(s)+\Si_x(s)^\top P(s)\Si_x(s)\\
\ns\ds\qq\qq\qq+Q(s)\Si_x(s)+\Si_x(s)^\top Q(s)+\BH_{xx}(s)\]ds+Q(s)
dW(s),\qq s\in[t,T],\\
\ns\ds P(T)=H_{xx}.\ea\right.\ee
Then
\bel{4.21}\ba{ll}
\ns\ds\dbE\(\h X_1^{i,\e,\bar\e}(T)^\top H_{xx}\h X_1^{i,\e,\bar\e}(T)+\int_t^T\h X_1^{i,\e,\bar\e}(r)^\top\BH_{xx}(r)\h X_1^{i,\e,\bar\e}(r)dr\)\\
\ns\ds=\dbE\int_t^T\Big\{\bar\xi_i^{\,\top}\[P(r)B_x(r)+B_x(r)^\top P(r)\]\bar\xi_i-\bar\xi_i^{\,\top}\Si_x(s)^\top P(r)\Si_x(r)\bar\xi_i+\si(\bar\t_i,r)^\top P(r)\si(\bar\t_i,r)\\
\ns\ds\qq-\bar\xi_i^{\,\top}\[Q(r)\Si_x(r)+\Si_x(r)^\top Q(r)\]\bar\xi_i+\bar\xi_i^{\,\top}Q(r)\Si_x(r)
\si(\bar\t_i,r)+\si(\bar\t_i,r)^\top\Si_x(r)^\top Q(r)\bar\xi_i\Big\}dr\ea\ee

\ep

Similar to BSDE \rf{BSDE1}, BSDE \rf{BSDE2} is also only depends on the optimal pair, and it is independent of $\xi^{i,\e,\bar\e}(\cd)$.

\ms

\it Proof. \rm The differential form of the equation for $\F^{i,\e,\bar\e}(\cd)=\h X_1^{i,\e,\bar\e}(\cd)\h X_1^{i,\e,\bar\e}(\cd)^\top$ is
\bel{XX*}\left\{\2n\ba{ll}
\ds d\F^{i,\e,\bar\e}(s)=\[B(s)\F^{i,\e,\bar\e}(s)+\F^{i,\e,\bar\e}(s)B(s)^\top
+\Si(s)\F^{i,\e,\bar\e}(s)\Si(s)^\top+\h b_{1,\bar\e}^i(s)\h X_1^{i,\e,\bar\e}(s)^\top\\
\ns\ds\qq\q+\h X_1^{i,\e,\bar\e}(s)\h b_{1,\bar\e}^i(s)^\top+\h\si_{1,\bar\e}^i(s)\h X_1^{i,\e,\bar\e}(s)^\top\Si(s)^\top
+\Si(s)\h X_1^{i,\e,\bar\e}(s)\h\si_{1,\bar\e}^i(s)^\top+\h\si^i_{1,\bar\e}(s)
\h\si^i_{1,\bar\e}(s)^\top\]ds\\
\ns\ds\qq\q+\[\Si(s)\F^{i,\e,\bar\e}(s)+\F^{i,\e,\bar\e}(s)\Si(s)^\top
+\h\si_{1,\bar\e}^i(s)\h X_1^{i,\e,\bar\e}(s)^\top+\h X_1^{i,\e,\bar\e}(s)
\h\si_{1,\bar\e}^i(s)^\top\]dW(s),\\
\ns\ds\F^{i,\e,\bar\e}(t)=0.\ea\right.\ee
By It\^o's formula, we have (suppressing $i,\e,\bar\e$ in $\F^{i,\e,\bar\e}(\cd)$)
$$\ba{ll}
\ns\ds d\dbE\big\{\tr[P(s)\F(s)]\big\}=\dbE\tr\Big\{\[-\(P(s)B_x(s)+B_x(s)^\top P(s)+\Si_x(s)^\top P(s)\Si_x(s)+Q(s)\Si_x(s)+\Si_x(s)^\top Q(s)\\
\ns\ds\qq\qq+\BH_{xx}(s)\)ds+Q(s)dW(s)\]\F(s)+P(s)\(B_x(s)\F(s)+\F(s)B_x(s)^\top
+\Si_x(s)\F(s)\Si_x(s)^\top\\
\ns\ds\qq\qq+\h b^i_{1,\bar\e}(s)\h X_1^{1,\e,\bar\e}(s)^\top+\h X_1^{1,\e,\bar\e}(s)\h b_{1,\bar\e}^i(s)^\top+\h\si_{1,\bar\e}^i(s)
\h X_1^{i,\e,\bar\e}(s)^\top\Si_x(s)^\top\\
\ns\ds\qq\qq+\Si_x(s)\h X_1^{i,\e,\bar\e}(s)
\h\si_{1,\bar\e}^i(s)^\top+\h\si_{1,\bar\e}^i(s)
\h\si_{1,\bar\e}^i(s)^\top\)ds\\
\ns\ds\qq\qq+P(s)\(\Si_x(s)\F(s)+\F(s)\Si_x(s)^\top+\h\si_{1,\bar\e}^i(s)
\h X_1^{i,\e,\bar\e}(s)^\top+\h X_1^{i,\e,\bar\e}(s)\h\si^i_{1,\bar\e}(s)^\top\)dW(s)\\
\ns\ds\qq\qq+Q(s)\(\Si_x(s)\F(s)+\F(s)\Si_x(s)^\top+\h\si_{1,\bar\e}^i(s)
\h X_1^{1,\e,\bar\e}(s)^\top+\h X_1^{i,\e,\bar\e}(s)\h\si_{1,\bar\e}^i(s)\)ds\Big\}\ea$$
$$\ba{ll}
\ns\ds=\dbE\tr\[-\BH_{xx}(s)\F(s)+P(s)\(\h b^i_{1,\bar\e}(s)\h X_1^{1,\e,\bar\e}(s)^\top+\h X_1^{1,\e,\bar\e}(s)\h b_{1,\bar\e}^i(s)^\top+\h\si_{1,\bar\e}^i(s)
\h X_1^{i,\e,\bar\e}(s)^\top\Si_x(s)^\top\\
\ns\ds\qq\qq+\Si_x(s)\h X_1^{i,\e,\bar\e}(s)
\h\si_{1,\bar\e}^i(s)^\top+\h\si_{1,\bar\e}^i(s)
\h\si_{1,\bar\e}^i(s)^\top\)ds+Q(s)\(\h\si_{1,\bar\e}^i(s)
\h X_1^{1,\e,\bar\e}(s)^\top+\h X_1^{i,\e,\bar\e}(s)\h\si_{1,\bar\e}^i(s)\)ds\]\\
\ns\ds=\dbE\[-\h X_1^{i,\e,\bar\e}(s)^\top\BH_{xx}(s)\h X_1^{i,\e,\bar\e}(s)+2\h X_1^{i,\e,\bar\e}(s)^\top P(s)\h b^i_{1,\bar\e}(s)+2\h X_1^{i,\e,\bar\e}(s)^\top\Si_x(s)^\top P(s)\h\si_{1,\bar\e}^i(s)
\\
\ns\ds\qq\qq+\h\si_{1,\bar\e}^i(s)^\top P(s)\h\si_{1,\bar\e}^i(s)+2\h X_1^{i,\e,\bar\e}(s)^\top Q(s)\h\si_{1,\bar\e}^i(s)\)ds\].\ea$$
Note that
$$\ba{ll}
\ns\ds\h X_1^{i,\e,\bar\e}(s)^\top P(s)\h b_{1,\bar\e}^i(s)=[X_1^{i,\e,\bar\e}(s)+\bar\xi_i{\bf1}_{[\bar\t_i,\bar\t_i+\bar\e)}
(s)]^\top P(s)B_x(s)\bar\xi_i{\bf1}_{[\bar\t_i,\bar\t_i+\bar\e)}
(s)\\
\ns\ds\qq\qq\qq\qq\q~=\bar\xi_i^{\,\top}P(s)B_x(s)\bar\xi_i
{\bf1}_{[\bar\t_i,\bar\t_i+\bar\e)}(s).\ea$$
Similarly,
$$\ba{ll}
\ns\ds\h X_1^{i,\e,\bar\e}(s)^\top\Si_x(s)^\top P(s)\h\si_{1,\bar\e}^i(s)=-[X_1^{i,\e,\bar\e}(s)+\bar\xi_i
{\bf1}_{[\bar\t_i,\bar\t_i+\bar\e)}]^\top\Si_x(s)^\top P(s)[\Si_x(s)\bar\xi_i-\si(\bar\t_i,s)]{\bf1}_{[\bar\t_i,
\bar\t_i+\bar\e)}(s)\\
\ns\ds=-\bar\xi_i^\top\Si_x(s)^\top P(s)[\Si_x(s)\bar\xi_i-\si(\bar\t_i,s)]{\bf1}_{[\bar\t_i,
\bar\t_i+\bar\e)}(s)\\
\ns\ds=\[-\bar\xi_i^{\,\top}\Si_x(s)^\top P(s)\Si_x(s)\bar\xi_i+\bar\xi^{\,\top}\Si_x(s)P(s)\si(\bar\t_i,s)\]{\bf1}_{[\bar\t_i,
\bar\t_i+\bar\e)}(s),\ea$$
and
$$\ba{ll}
\ns\ds\h\si_{1,\bar\e}^i(s)^\top P(s)\h\si_{1,\bar\e}^i(s)=[\Si_x(s)\bar\xi_i-
\si(\bar\t_i,s)]^\top P(s)[\Si_x(s)\bar\xi_i-\si(\bar\t_i,s)]{\bf1}_{[\bar\t_i,
\bar\t_i+\bar\e)}(s)\\
\ns\ds=\[\bar\xi_i^{\,\top}\Si_x(s)^\top P(s)\Si_x(s)\bar\xi_i-\si(\bar\t_i,s)^\top P(s)\Si_x(s)\bar\xi_i-\bar\xi_i^{\,\top}\Si_x(s)^\top P(s)\si(\bar\t_i,s)+\si(\bar\t_i,s)^\top P(s)\si(\bar\t_i,s)\]{\bf1}_{[\bar\t_i,\bar\t_i+\bar\e)}(s),\\
\ns\ds\h X_1^{i,\e,\bar\e}(s)^\top Q(s)\h \si_{1,\bar\e}^i(s)=-[X_1^{i,\e,\bar\e}(s)+\bar\xi_i{\bf1}_{[\bar\t_i,\bar\t_i+\bar\e)}
(s)]^\top Q(s)[\Si_x(s)\bar\xi_i-\si(\bar\t_i,s)]{\bf1}_{[\bar\t_i,\bar\t_i+\bar\e)}
(s)\\
\ns\ds\qq=\[-\bar\xi_i^{\,\top}Q(s)\Si_x(s)\bar\xi_i+\bar\xi_i^{\,\top}Q(s)\Si_x(s)
\si(\bar\t_i,s)\]{\bf1}_{[\bar\t_i,\bar\t_i+\bar\e)}
(s).\ea$$
Hence,
$$\ba{ll}
\ns\ds\dbE\[2\h X_1^{i,\e,\bar\e}(s)^\top P(s)\h b^i_{1,\bar\e}(s)+2\h X_1^{i,\e,\bar\e}(s)^\top\Si_x(s)^\top P(s)\h\si_{1,\bar\e}^i(s)
\\
\ns\ds\qq\qq+\h\si_{1,\bar\e}^i(s)^\top P(s)\h\si_{1,\bar\e}^i(s)+2\h X_1^{i,\e,\bar\e}(s)^\top Q(s)\h\si_{1,\bar\e}^i(s)\]\\
\ns\ds=\dbE\Big\{2\bar\xi_i^{\,\top}P(s)B_x(s)\bar\xi_i+2\[-\bar\xi_i^{\,\top}\Si_x(s)^\top P(s)\Si_x(s)\bar\xi_i+\bar\xi^{\,\top}\Si_x(s)P(s)\si(\bar\t_i,s)\]\\
\ns\ds\qq+\[\bar\xi_i^{\,\top}\Si_x(s)P(s)\Si_x(s)\bar\xi_i-\si(\bar\t_i,s)^\top P(s)\Si_x(s)\bar\xi_i-\bar\xi_i^{\,\top}\Si_x(s)^\top P(s)\si(\bar\t_i,s)+\si(\bar\t_i,s)^\top P(s)\si(\bar\t_i,s)\]\\
\ns\ds\qq+2\[-\bar\xi_i^{\,\top}Q(s)\Si_x(s)\bar\xi_i+\bar\xi_i^{\,\top}Q(s)\Si_x(s)
\si(\bar\t_i,s)\]\Big\}{\bf1}_{\bar\t_i,\bar\t_i+\bar\e)}(s)\\
\ns\ds=\dbE\Big\{\bar\xi_i^{\,\top}\[P(s)B_x(s)+B_x(s)^\top P(s)\]\bar\xi_i-\bar\xi_i^{\,\top}\Si_x(s)^\top P(s)\Si_x(s)\bar\xi_i+\si(\bar\t_i,s)^\top P(s)\si(\bar\t_i,s)\\
\ns\ds\qq-\bar\xi_i^{\,\top}\[Q(s)\Si_x(s)+\Si_x(s)^\top Q(s)\]\bar\xi_i+\bar\xi_i^{\,\top}Q(s)\Si_x(s)
\si(\bar\t_i,s)+\si(\bar\t_i,s)^\top\Si_x(s)^\top Q(s)\bar\xi_i\Big\}{\bf1}_{\bar\t_i,\bar\t_i+\bar\e)}(s)\ea$$
Thus, integrating the above over $[t,T]$ to get \rf{4.21}. \endpf

\ms

From the above, the variational inequality reads
\bel{ineq4}\ba{ll}
\ns\ds0\les J(\t_0,t,x;\xi^{i,\e,\bar\e}(\cd))-J(\t_0,t,x;\bar\xi(\cd))\\
\ns\ds\q=\dbE\Big\{\int_t^T\[\bar\xi_i^{\,\top}\(P(r)B_x(r)+B_x(r)^\top P(r)-\Si_x(\t)^\top P(r)\Si_x(r)-Q(r)\Si_x(r)-\Si_x(r)^\top Q(r)\)\bar\xi_i\\
\ns\ds\qq\qq+\si(\bar\t_i,r)^\top P(r)\si(\bar\t_i,r)+\bar\xi_i^{\,\top}Q(r)\Si_x(r)
\si(\bar\t_i,r)+\si(\bar\t_i,r)^\top\Si_x(r)^\top Q(r)\bar\xi_i\\
\ns\ds\qq\qq-\lan\BH_x(r),\bar\xi_i\ran-\BH(r)\]
{\bf1}_{[\bar\t_i,\bar\t_i+\bar\e)}(r)dr\\
\ns\ds\qq\qq+\bar\e\(\ell_\t(\bar\t_i,\bar\xi_i)+H_x\z^i(T)+\int_t^T[\lan\BH_x(r),\z^i(r)\ran+
g^i_\t(r)]dr\)\\
\ns\ds\qq\qq+\e\(H_x+\ell_\xi(\bar\t_i,\bar\xi_i)+\int_t^TG_x(r)
{\bf1}_{[\bar\t_i,T]}(r)dr\)\xi_i\Big\}+o(\e+\bar\e)\ea\ee
When we replace $\xi^{i,\bar\e,\e}(\cd)$ by $\xi^{i,\e,-\bar\e}(\cd)$ (see \rf{xi^e} and \rf{xi-}), the variational inequality becomes
\bel{ineq5}\ba{ll}
\ns\ds0\les J(\t_0,t,x;\xi^{i,\e,-\bar\e}(\cd))-J(\t_0,t,x;\bar\xi(\cd))\\
\ns\ds\q=\dbE\Big\{\int_t^T\[\bar\xi_i^{\,\top}\(P(r)B_x(r)+B_x(r)^\top P(r)-\Si_x(\t)^\top P(r)\Si_x(r)-Q(r)\Si_x(r)-\Si_x(r)^\top Q(r)\)\bar\xi_i\\
\ns\ds\qq\qq+\si(\bar\t_i,r)^\top P(r)\si(\bar\t_i,r)+\bar\xi_i^{\,\top}Q(r)\Si_x(r)
\si(\bar\t_i,r)+\si(\bar\t_i,r)^\top\Si_x(r)^\top Q(r)\bar\xi_i\\
\ns\ds\qq\qq-\lan\BH_x(r),\bar\xi_i\ran-\BH(r)\]
{\bf1}_{[\bar\t_i,\bar\t_i+\bar\e)}(r)dr\\
\ns\ds\qq\qq-\bar\e\(\ell_\t(\bar\t_i,\bar\xi_i)+H_x\z^i(T)+\int_t^T[\lan\BH_x(r),\z^i(r)\ran+
g^i_\t(r)]dr\)\\
\ns\ds\qq\qq+\e\(H_x+\ell_\xi(\bar\t_i,\bar\xi_i)+\int_t^TG_x(r)
{\bf1}_{[\bar\t_i,T]}(r)dr\)\xi_i\Big\}+o(\e+\bar\e)\ea\ee
Hence, we have the following Pontryagin type maximum principle for our Problem (IC).

\bt{MP} \sl Let {\rm(H1)--(H2)} and {\rm(H4)} hold. Let $(\t,t,x)\in\D_*[0,T]\times\dbR^n$ be given and $\bar\xi(\cd)\in\BBXi[t,T]$ be an optimal impulse control of from \rf{bar xi}. Let $(Y(\cd),Z(\cd))$ and $(P(\cd), Q(\cd))$ be the adapted solutions to BSDEs \rf{BSDE1} and \rf{BSDE2} respectively. Then, for $i\in\{2,\cds,\k(\bar\xi(\cd))-1\}$, or $\bar\t_1>t$, $\bar\t_{\k(\bar\xi(\cd))}<T$, $\xi_i=\eta-\bar\xi_i$ with $\eta\in\sK$, almost surely, it holds:
\bel{MP1}\ba{ll}
\ns\ds\q\bar\xi_i^{\,\top}\[P(\bar\t_i)B_x(\bar\t_i)+B_x(\bar\t_i)^\top P(\bar\t_i)-\Si_x(\bar\t_i)^\top P(\bar\t_i)\Si_x(\bar\t_i)-Q(\bar\t_i)\Si_x(\bar\t_i)-\Si_x(\bar\t_i)^\top Q(\bar\t_i)\\
\ns\ds\qq\qq+Q(\bar\t_i)\Si_x(\bar\t_i)\si(\bar\t_i,\bar\t_i)
+\si(\bar\t_i,\bar\t_i)^\top\Si_x(\bar\t_i)^\top Q(\bar\t_i)\]\bar\xi_i\\
\ns\ds\qq\qq+\si(\bar\t_i,\bar\t_i)^\top P(\bar\t_i)\si(\bar\t_i,\bar\t_i)-\lan\BH_x(\bar\t_i),\bar\xi_i\ran
-\BH(\bar\t_i)\ges0,\ea\ee
\bel{MP2}\(H_x+\ell_\xi(\bar\t_i,\bar\xi_i)+\int_t^TG_x(r)
{\bf1}_{[\bar\t_i,T]}(r)dr\)\xi_i\ges0,\ee
\bel{MP3}\ell_\t(\bar\t_i,\bar\xi_i)+H_x\z^i(T)+\int_t^T[\lan\BH_x(r),\z^i(r)\ran+
g^i_\t(r)]dr=0.\ee
In the case $\bar\t_1=t$, \rf{MP1} and \rf{MP2} remain true, and \rf{MP3} becomes
\bel{MP4}\ell_\t(\bar\t_i,\bar\xi_i)+H_x\z^i(T)+\int_t^T[\lan\BH_x(r),\z^i(r)\ran+
g^i_\t(r)]dr\ges0.\ee
In the case $\bar\t_{\k(\bar\xi(\cd))}=T$, \rf{MP1} and \rf{MP2} remain true, and \rf{MP3} becomes
\bel{MP5}\ell_\t(\bar\t_i,\bar\xi_i)+H_x\z^i(T)+\int_t^T[\lan\BH_x(r),\z^i(r)\ran+
g^i_\t(r)]dr\les0.\ee

\et

Note that if we only (convexly) perturb the sizes $\bar\xi_i$ of the optimal impulse control, one only gets inequality \rf{MP2}, and the second adjoint equation \rf{BSDE2} is not necessary.

\section{Extensions and Conclusions}

In this section, we first present some possible extensions.

\ms

Consider the system (state and cost functional) also contains
other types of controls. Thus, our controlled state equation looks like
\bel{state'}\ba{ll}
\ns\ds X(s)=\[x+\int_t^sb(\t_0,\t,X(\t),u(\t))d\t+\int_t^s\si(\t_0,\t,X(\t),u(\t))
dW(\t)\\
\ns\ds\qq\qq+\sum_{j\ges1}\(\int_{\t_j}^sb(\t_j,\t,X(\t),u(\t))d\t
+\int_{\t_j}^s\si(\t_j,\t,X(\t),u(\t))dW(\t)\)\]{\bf1}_{[\t_j,T]}(s)\\
\ns\ds\qq\qq\qq\qq\qq\qq\qq\qq\qq\qq\qq+\xi(s),\q s\in[t,T],\ea\ee
with recursive functional
\bel{J**}J(\t_0,t,x;\xi(\cd))=Y(\t_0,t,t,\xi(\cd)),\ee
where $(Y(\cd),Z(\cd\,,\cd))$ is the adapted solution to the following BSVIE
\bel{BSVIE*}\ba{ll}
\ns\ds Y(s)=h(X(T))+\sum_{j\ges1}\ell(\t_j,X(\t_j-0),\xi_j){\bf1}_{[\t_j,T]}(s)\\
\ns\ds\qq\qq+\int_s^T\(g(\t_0,\t,X(\t),u(\t))+\sum_{j\ges1}
g(\t_j,\t,X(\t),u(\t))\)d\t-\int_s^TZ(s,\t)dW(\t),\q s\in[t,T],\ea\ee
where $g:\D_*[0,T]\times\dbR^n\times U\to\dbR_+$, $h:\dbR^n\to\dbR_+$, and $\ell:
[0,T]\times\dbR^n\times\sK\to\dbR_+$. Then, using the arguments before, we could derive the continuity of the value function, the corresponding HJB equation which should be a quasi-variational inequality. Further, one is able to show that the value function is the unique viscosity solution of its HJB equation.
Next, by making perturbations of the sizes and impulse moments of the optimal impulse control, together with general perturbation of the optimal (continuously active) control, one could obtain the corresponding maximum principle. We leave the details (maybe complicated) to interested readers.

\ms

It is possible that switching control is also get involved. In that case, the
results of perturbation with respect to the optimal switching moments are unknown. But, we suspect that the technique introduced in the current paper might be helpful. We encourage interested readers can complete this task.

\ms

We also point out that the system is allowed to depend on a Markov chain governed regime-switching. Then more interesting phenomena might appear. We expect the forms of corresponding quasi-variational inequality, and its maximum principle look like.

\ms

To conclude this paper, we would like to mention two major contributions. The first, we are able to provide a framework so that the impulse control problem is allowed to have a changing running cost. The second is to perturb the optimal impulse moments so that we derive a complete version of the maximum principle.

\end{document}